\documentclass[preprint,aihp]{imsart}

\usepackage[]{graphicx} 
\usepackage{amsfonts,float}

\usepackage{amssymb}

\DeclareSymbolFont{sfletters}{OML}{cmbrm}{m}{it}

\DeclareMathSymbol{\salpha}{\mathord}{sfletters}{"0B}
\usepackage{mathrsfs}
\RequirePackage{natbib}
\RequirePackage[colorlinks,citecolor=blue,urlcolor=blue]{hyperref}
\usepackage{overpic}
\makeatletter
\newcommand{\oset}[3][1ex]{%
  \mathrel{\mathop{#3}\limits^{
    \vbox to#1{\kern-2\ex@
    \hbox{$\scriptstyle#2$}\vss}}}}
\makeatother

\usepackage{mathtools}
\RequirePackage{natbib}
\usepackage{comment}
\usepackage{enumerate}

\usepackage{eufrak}

\usepackage{array,varwidth}

\usepackage{overpic}
\usepackage[framemethod=tikz,innertopmargin=7pt]{mdframed}
\usepackage{graphicx} 
\usepackage{subcaption}
\usepackage{MnSymbol}
\usepackage{placeins} 
\usepackage{scalerel}

\usepackage[cal=zapfc]{mathalfa}

\makeatletter

\def \@cflt{%
    \let \@elt \@comflelt
    \setbox\@tempboxa \vbox{}%
    \@toplist
\setbox\z@\vsplit\@outputbox to 0.5\ht\@outputbox
    \setbox\@outputbox \vbox{%
                             \boxmaxdepth \maxdepth
                             \unvbox\z@
                             \vskip .5\textfloatsep
                             \unvbox\@tempboxa
                             \vskip -\floatsep
                             \topfigrule
                             \vskip .5\textfloatsep
                             \unvbox\@outputbox
                             }%
    \let\@elt\relax
    \xdef\@freelist{\@freelist\@toplist}%
    \global\let\@toplist\@empty
}

\makeatother

\usepackage{booktabs,longtable,tabu} 
\usepackage{multirow} 
\usepackage{float}
\usepackage{makecell}
\usepackage{dcolumn}

\newcommand*\widebar[1]{%
   \vbox{%
     \hrule height 0.5pt
     \kern0.25ex
     \hbox{%
       \kern-0.15em
       \ifmmode#1\else\ensuremath{#1}\fi
       \kern-0.15em
     }
   }
}

\RequirePackage[OT1]{fontenc}
\RequirePackage{amsthm,amsmath}

\startlocaldefs
\numberwithin{equation}{section}
\theoremstyle{plain}
\newtheorem{theorem}{Theorem}[section]
\newtheorem{lemma}{Lemma}[section]
\newtheorem{proposition}{Proposition}[section]
\newtheorem{assumption}{Assumption}[section]

\usepackage{algorithm}
\usepackage{algorithmic}

\theoremstyle{remark}

\newcommand{\R}{\mathbb{R}}

\renewcommand{\hat}[1]{\widehat{#1}}

\newcommand{\vunderline}[1]{\underline{#1\mkern-2mu}\mkern2mu }

\newcommand{\ok}{\boldsymbol{\pmb{\otimes}}}

\newcommand{\I}{\text{\textsc{I}}}
\newcommand{\II}{\text{\textsc{II}}}

\newcommand{\III}{\text{\textsc{III}}}

\newcommand{\dK}{d_{\textup{K}}}

\newcommand{\J}{{\mathsf{J}}}

\renewcommand{\P}{\textup{\textbf{P}}}
\newcommand{\E}{\textup{\textbf{E}}}

\newcommand{\LL}{\mathcal{L}}

\newcommand{\e}{\epsilon}

\newcommand{\cov}{\textup{cov}}
\newcommand{\tr}{\operatorname{tr}}

\newcommand{\var}{\operatorname{var}}

\newcommand{\ttop}{^{\top}}

\newcommand{\op}{_{\textup{op}}}
\renewcommand{\dim}{{\textup{dim}}}

\newcommand{\ts}{\textstyle}

\newcommand{\diag}{\textsf{Diag}}

\endlocaldefs

\usepackage[T1]{fontenc}
\begin{document}

\begin{frontmatter}
\title{ Improved Rates of Bootstrap Approximation for the Operator Norm: A Coordinate-Free Approach}
\runtitle{Improved Rates of Bootstrap Approximation for the Operator Norm}

\begin{aug}
\author[A]{\fnms{Miles E.}~\snm{Lopes}\ead[label=e1]{melopes@ucdavis.edu}},

\runauthor{M.~E. Lopes}

%
\address[A]{University of California, Davis\printead[presep={,\ }]{e1}}

\end{aug}

\begin{abstract}
Let $\hat\Sigma=\frac{1}{n}\sum_{i=1}^n X_i\otimes X_i$ denote the sample covariance operator of  centered i.i.d.~observations $X_1,\dots,X_n$ in a real separable Hilbert space, and let $\Sigma=\E(X_1\otimes X_1)$. The focus of this paper is to understand how well the bootstrap can approximate the distribution of the operator norm error $\sqrt n\|\hat\Sigma-\Sigma\|\op$, in settings where the eigenvalues of $\Sigma$ decay as $\lambda_j(\Sigma)\asymp j^{-2\beta}$ for some fixed parameter $\beta>1/2$. Our main result shows that the bootstrap can approximate the distribution of $\sqrt n\|\hat\Sigma-\Sigma\|\op$ at a rate of order  $n^{-\frac{\beta-1/2}{2\beta+4+\e}}$ with respect to the Kolmogorov metric, for any fixed $\e>0$. In particular, this shows that the bootstrap can achieve near $n^{-1/2}$ rates in the regime of large $\beta$---which substantially improves on previous near $n^{-1/6}$ rates in the same regime. In addition to obtaining faster rates, our analysis leverages a fundamentally different perspective based on coordinate-free techniques.  Moreover, our result holds in greater generality, and we propose a model that is compatible with both elliptical and Mar\v{c}enko-Pastur models in high-dimensional Euclidean spaces, which may be of independent interest.
 \end{abstract}


\end{frontmatter}

\section{Introduction}\label{sec:intro}
Let $X_1,\dots,X_n$ be centered i.i.d.~observations in a real separable Hilbert space $(\mathbb{H},\langle \cdot,\cdot\rangle)$, and let the sample and population covariance operators be denoted respectively as \smash{$\hat\Sigma=\frac{1}{n}\sum_{i=1}^n X_i\otimes X_i$} and $\Sigma=\E(X_1\otimes X_1)$. In the fields of high-dimensional statistics and functional data analysis, many strands of research are connected to the study of the operator norm error, denoted
\begin{equation*}
T_n = \sqrt n \|\hat\Sigma-\Sigma\|\op.
\end{equation*}
Indeed, the operator norm is a standard measure of error in the extensive literature on covariance estimation \citep[e.g.][]{Rudelson:1999,Bickel2:2008,Cai:2010,Adamczak:2011,Lounici:2014, Bunea:Xiao:2015,Koltchinskii:Bernoulli,Minsker:2017,Zhivotovskiy}. 
Furthermore, the statistic $T_n$ plays an important role in principal components analysis, where it appears frequently in various error bounds for sample eigenvalues and eigenvectors
(e.g.~\citep[][\textsection 4.2]{Bosq},~\citep{Samworth:2014}). More broadly still, the statistic $T_n$ is closely related to errors that arise in the analysis of randomized sketching algorithms for large-scale matrix computations~\citep[e.g.][]{Tropp:Acta,Lopes:Bernoulli,RFF}.

 Despite the prevalence of the statistic $T_n$, the standard tools for approximating its distribution have some essential limitations. 
At a high level, the bulk of existing results describing $T_n$ can be organized into two types. In the direction of high-dimensional statistics, the large body of work on covariance estimation has developed non-asymptotic tail bounds for $T_n$, as cited earlier. Alternatively, in the direction of functional data analysis, the relevant theory has focused mainly on the limiting distribution for the operator $\sqrt n (\hat\Sigma-\Sigma)$ as $n\to\infty$, which implicitly describes $T_n$ via the continuous mapping theorem~\citep[e.g.][]{Dauxois:1982,Bosq,Mas:2002,Mas:2006,Hsing:2015,Paparoditis:2016}. 

While both of these types of results yield valuable insights, they do not provide a complete understanding of the fluctuations of $T_n$. For instance, asymptotic distributional approximations generally do not provide much guidance about how quickly accuracy improves as $n\to\infty$. Meanwhile, non-asymptotic tail bounds are often impractical for approximating the numerical values of tail probabilities, because such bounds are usually stated in terms of unspecified or conservative constants. On this point, it is worth clarifying that non-asymptotic distributional approximations often also involve unspecified constants, but they are less of a practical concern. To illustrate the difference, a tail bound for $T_n$ typically has the form $\P(T_n> t)\leq  C b(t,n,\Sigma)$, where $b$ is a function and $C$ is a constant. By contrast, non-asymptotic distributional approximation results often show that an event of the form $\sup_{t\in\R} |\P(T_n\leq t)-\hat F_n(t)|\leq C b(n,\Sigma)$ holds with high probability, where $\hat F_n$ is an estimate for the distribution function of $T_n$. Hence, in the latter context, the constant $C$ plays an essentially theoretical role since it does not interfere with the practical use of $\hat F_n(t)$, whereas the tail bound $Cb(t,n,\Sigma)$ cannot be applied numerically without specifying a particular value of $C$.

In light of these issues, there is substantial interest in establishing \emph{rates} of distributional approximation. Recently, this interest has stimulated a line of work on rates of bootstrap and Gaussian approximations for the distributions of $T_n$ and other statistics related to $\hat\Sigma$~\citep{KoltchinskiiAOS:2017,Zhou:2018,Lopes:Bernoulli,Naumov:2019,Silin,Yao:Lopes,jirakboot}. Here, we only discuss the papers~\citep{Zhou:2018,Lopes:Bernoulli} from this list, since they deal specifically with $T_n$, whereas the others deal with qualitatively different statistics arising from $\hat\Sigma$. When the data are $p$-dimensional Euclidean vectors, i.e.~$\mathbb{H}=\R^p$, the papers ~\citep{Zhou:2018,Lopes:Bernoulli} establish non-asymptotic high-probability upper bounds on $\sup_{t\in\R}|\P(T_n\leq t)-\P(T_n^*\leq t|X)|$, where $T_n^*$ is a bootstrap version of $T_n$, and $\P(\cdot\,|X)$ refers to probability that is conditional on $X_1,\dots,X_n$. In a setting where the eigenvalues of $\Sigma$ are not assumed to have any decay structure, the paper~\citep[Theorem 2.3]{Zhou:2018} obtains a bound of order $p^{9/8}/n^{1/8}$, up to factors of $\log(n)$. On the other hand, in the distinct setting where the eigenvalues of $\Sigma$ decay according to $\lambda_j(\Sigma)\asymp j^{-2\beta}$ for some fixed parameter $\beta>1/2$,  the paper~\citep[Theorem 2.1]{Lopes:Bernoulli} obtains a bound of order $n^{-\frac{\beta-1/2}{6\beta+4}}$, up to factors of $\log(n)$.

 Although the result in~\citep{Lopes:Bernoulli} has the virtue of being applicable to very high-dimensional data with $p\gg n$, it still has several key drawbacks---which the current paper seeks to overcome. Most notably, the rate $n^{-\frac{\beta-1/2}{6\beta+4}}$ can \emph{never be faster than $n^{-1/6}$}, for any value of $\beta>1/2$.  At a heuristic level, such a rate seems to leave substantial room for improvement, because if $\beta$ is large, then the data will be concentrated near a low-dimensional subspace of $\mathbb{H}$, which should seemingly enable the bootstrap to achieve rates of approximation that are near $n^{-1/2}$. Another drawback is that the result in~\citep{Lopes:Bernoulli} requires the data to take values in finite-dimensional Euclidean spaces, which is incongruous with the fact that the result is dimension-free. In particular, the result cannot be applied to standard types of functional data that reside in infinite-dimensional Hilbert spaces.\\

\noindent\textbf{Contributions.} In relation to the points just discussed, the main contributions of this paper are briefly summarized as follows:
\begin{enumerate}
\item Under conditions more general than those in~\citep{Lopes:Bernoulli}, our main result gives a high-probability upper bound on $\sup_{t\in\R}|\P(T_n\leq t)-\P(T_n^*\leq t|X)|$ that is of order $n^{-\frac{\beta-1/2}{2\beta+4+\e}}$ for any fixed $\e>0$. In particular, this rate is near $n^{-1/2}$ in the regime of large $\beta$, and in fact,  it is faster than $n^{-\frac{\beta-1/2}{6\beta+4}}$ for every value of $\beta>1/2$.\\[-0.2cm]

\item We allow the observations $X_1,\dots,X_n$ to take values in an abstract separable Hilbert space $\mathbb{H}$. Our setup also allows $\mathbb{H}$ to vary with $n$, which makes it possible to handle both functional data and high-dimensional data, with $\mathbb{H}=\R^p$ and $p\gg n$, in a unified way.\\[-0.2cm]

\item We develop a data-generating model that is compatible with both Mar\v{c}enko-Pastur models and elliptical models in high dimensions (see Propositions~\ref{prop:elliptical} and~\ref{prop:MP}). This model has the favorable property that it can be defined by a concise set of assumptions within the setting of an abstract separable Hilbert space. Likewise, this model may be of independent interest in both high-dimensional statistics and functional data analysis.\\[-0.2cm]

\item Whereas the papers~\citep{Lopes:Bernoulli} and~\citep{Zhou:2018} represent  $T_n$ as $\sup_{\|u\|=1}|\frac{1}{\sqrt n}\sum_{i=1}^n \langle u,X_i\rangle^2 - \E \langle u,X_i\rangle^2|$ and analyze it as a supremum of an empirical process, 
our work takes a fundamentally different perspective based on \emph{coordinate-free techniques}, as explained below.
\end{enumerate}

\noindent\textbf{Proof techniques.} Our approach begins with a preliminary dimension-reduction step, based on an approximation of the form $ \|\hat\Sigma-\Sigma\|\op \approx \|\frac{1}{n}\sum_{i=1}^n Y_iY_i\ttop -\E(Y_iY_i\ttop)\|\op$, where $Y_1,\dots,Y_n$ are random vectors in $\R^{k}$ and $k$ may increase with $n$. Once this has been done, the first key idea is to view $Y_1Y_1\ttop,\dots,Y_nY_n\ttop$ as \emph{vectors} of dimension \smash{$\frac{1}{2}k(k+1)$,} rather than as symmetric matrices of size $k\times k$. The immediate benefit of this viewpoint is that it enables us to leverage well-developed Berry-Esseen bounds that can handle general norms of sums of Euclidean vectors~\citep[e.g.][]{Fang:Koike:Balls,Bentkus:2003}.

However, the price to pay for this benefit is that we need to deal with the covariance matrix of the ``vector'' $\frac{1}{n}\sum_{i=1}^n Y_iY_i\ttop -\E(Y_iY_i\ttop)$. Here, it is crucial to notice that when this sum is viewed as a vector in $\R^{k(k+1)/2}$, its covariance matrix is a complicated object. Indeed, if we let $\mathsf{K}$ denote the covariance matrix of $VV\ttop$ for a generic random vector $V\in\R^k$, then the matrix  $\mathsf{K}$ is of size \smash{$\frac{1}{2}k(k+1)\times \frac{1}{2}k(k+1)$} and has entries that are indexed by 4-tuples $(j_1,j_2,j_3,j_4)\in\{1,\dots,k\}^4$. To make matters even more challenging, we will sometimes need to work with the inverse $\mathsf{K}^{-1}$, which is prohibitive when using an entrywise representation of $\mathsf{K}$.

The second key idea is to recognize the advantages of coordinate-free techniques. That is, rather than thinking of $\mathsf{K}$ as a matrix, we think of it as an abstract linear operator. 
 To explain this more precisely, first let $\mathbb{S}^{k\times k}$ denote the Hilbert space  of symmetric matrices of size $k\times k$, equipped with the inner product $\llangle A, B\rrangle=\tr(AB)$. Then, we can define $\mathsf{K}$ in a coordinate-free manner as the linear operator from $\mathbb{S}^{k\times k}$ to $\mathbb{S}^{k\times k}$ that satisfies 
\begin{equation}\label{eqn:abcov}
\llangle A, \mathsf{K}(B)\rrangle \ = \ \cov(\llangle A, VV\ttop \rrangle, \llangle B, VV\ttop \rrangle)
\end{equation}
 for all $A,B\in\mathbb{S}^{k\times k}$, where $\cov(\cdot,\cdot)$ is the usual covariance for scalar random variables~\citep[][]{Eaton}. 
 The main advantage of this viewpoint is that it avoids the complexity of $\mathsf{K}$ at an entrywise level. In addition, there are two other further advantages. First, it turns out that much of the information we need to extract from the inverse $\mathsf{K}^{-1}$ can be obtained by minimizing and maximizing the function $A\mapsto \llangle A, \mathsf{K}(A)\rrangle $ over certain subsets of $\mathbb{S}^{k\times k}$. Second, the inner product $\llangle A, \mathsf{K}(A)\rrangle $ is often amenable to direct calculation, since it is the same as the variance of a quadratic form, $\var(V\ttop A V)$. This being said, we should also clarify that our use of coordinate-free techniques is more utilitarian than aesthetic, and we do not hesitate to use coordinates in other aspects of our work.
 
One more item to expand on is the dimension-reduction step mentioned earlier. This involves projecting $\mathbb{H}$-valued data into $k$-dimensional Euclidean spaces with $k\to\infty$ as $n\to\infty$, which is an essential technique in many aspects of functional data analysis~\citep[e.g.][]{ramsaysilverman,Bosq,CaiHall,Horvath,HorvathJMVA}. In our work, this technique allows us to analyze $\|\frac{1}{n}\sum_{i=1}^n Y_iY_i\ttop -\E(Y_iY_i\ttop)\|\op$ as a convenient surrogate for $\|\hat\Sigma-\Sigma\|\op$ (see Lemma~\ref{lem:covcouple}).
  Perhaps the main point of contrast between our use this technique and more common uses in functional data analysis is that we focus on the operator norm. This is because the Frobenius (Hilbert-Schmidt) norm tends to play a more dominant role in the theory of functional data. Another notable feature of our work is that we account for the projected dimension $k$ in a quantitative manner, which is made possible by results from high-dimensional probability. In this respect, our main tool is Lemma~\ref{lem:op} (established in~\citep[][Theorem 2.2]{Lopes:Bernoulli}), which provides a bound on the $q$th moment $\E\|\hat\Sigma-\Sigma\|\op^q$ in the general setting of separable Hilbert spaces and allows $q\to\infty$ as $n\to\infty$. At the same time, it should be emphasized that this bound on $\E\|\hat\Sigma-\Sigma\|\op^q$ is but one of many related results~\citep[e.g.][]{mendelson2005,mendelson2006,oliveira,troppuser,Koltchinskii:Bernoulli} that developed from the seminal work~\citep{Rudelson:1999} based on non-commutative Khintchine inequalities.\\


\noindent\textbf{Notation and terminology.} The Hilbert norm on $\mathbb{H}$ is denoted $\|\cdot\|$, and the direct sum of $\mathbb{H}$ with itself is denoted $\mathbb{H}\oplus\mathbb{H}$. The tensor product $x\otimes y$ of any two vectors $x,y\in \mathbb{H}$ is defined as the linear operator that satisfies $(x\otimes y)z = \langle y,z\rangle x$ for all $z\in\mathbb{H}$. All Euclidean vectors are regarded as column vectors, except when they are appended with a $\top$. 
When $x$ and $y$ are Euclidean vectors, we will use $x\otimes y$ and $xy\ttop$ interchangeably, depending on notational convenience. As an important clarification, the expression $A\otimes B$ does \emph{not} denote the Kronecker product of $A,B\in\mathbb{S}^{k\times k}$. Rather, we use $A\otimes B$ to denote the tensor product of $A$ and $B$, viewed as elements of the Hilbert space $\mathbb{S}^{k\times k}$. When referring to the Kronecker product of $A$ and $B$, we will write $A\ok B$, which is the linear operator on $\mathbb{S}^{k\times k}$ that satisfies $(A\ok B)C= A C B$ for all $C\in\mathbb{S}^{k\times k}$. Additional background on coordinate-free definitions of tensor and Kronecker products in the context of multivariate statistics can be found in~\citep{Eaton}.

If $\mathsf{L}:\mathbb{H}\to\mathbb{H}$  is a linear operator, then its operator norm is defined by $\|\mathsf{L}\|\op=\sup_{\|u\|=1}\|\mathsf{L}u\|$, and its Frobenius norm is defined by $\|\mathsf{L}\|_F^2=\sum_{j}\|\mathsf{L}\phi_j\|^2$, where $\{\phi_j\}$ is any orthonormal basis for $\mathbb{H}$. If $\mathsf{L}$ is compact and positive-semidefinite, or if $\mathsf{L}$ is self-adjoint and $\dim(\mathbb{H})<\infty$, then we write $\lambda_1(\mathsf{L})\geq \lambda_2(\mathsf{L})\geq \cdots$ for the sorted eigenvalues, with $\lambda_{\max}(\mathsf{L})$ and $\lambda_{\min}(\mathsf{L})$ denoting the maximum and minimum eigenvalues in the latter case. In the Euclidean space $\R^k$, the standard basis vectors are $e_1,\dots,e_k$, and the ordinary Euclidean norm is $\|\cdot\|_2$. For real numbers $a_1,\dots,a_k$, we define $\diag(a_1,\dots,a_k)$ as the $k\times k$ diagonal matrix whose $j$th diagonal entry is $a_j$. For a matrix $A\in\mathbb{S}^{k\times k}$, we define $\diag(A)$ as the $k\times k$ diagonal matrix obtained by setting the off-diagonal entries of $A$ to 0. The identity matrix is denoted by $I$, with its size being understood from context.

If $U$ is a scalar random variable, then its $L^q$ norm is defined by $\|U\|_{L^q}=(\E|U|^q)^{1/q}$ for any $q\geq 1$. The basic inequality $\P(|U|\geq e\|U\|_{L^q})\leq e^{-q}$ will often be referred to as Chebyshev's inequality. The expression $\mathcal{L}(U)$ refers to the distribution of $U$. If $V$ is another scalar random variable, then $U\overset{{\scriptsize\mathcal{L}}}{=} V$ means that $U$ and $V$ are equal in distribution. As a more compact notation for the Kolmogorov metric, we write $d_{\textup{K}}(\mathcal{L}(U),\mathcal{L}(V))=\sup_{t\in\R}|\P(U\leq t)-\P(V\leq t)|$. A random vector is said to be isotropic if its covariance matrix is equal to the identity matrix.

If $\{a_n\}$ and $\{b_n\}$ are sequences of non-negative real numbers, the relation $a_n\lesssim b_n$ means that there is a constant $c>0$ not depending on $n$ such that $a_n\leq c \,b_n$ holds for all large $n$. If both of the relations $a_n\lesssim b_n $ and $b_n\lesssim a_n$ hold, then we write $a_n\asymp b_n$. When introducing constants that do not depend on $n$, we will often re-use the symbol $c$, so that its value may change at each appearance. For the maximum and minimum of a pair of real numbers $a$ and $b$, we write $a\vee b=\max\{a,b\}$ and $a\wedge b=\min\{a,b\}$.

\section{Main result}

Our work will be done with a sequence of models implicitly embedded in a triangular array, whose rows are indexed by $n$. The space $\mathbb{H}$ and all model parameters are allowed to vary freely with $n$, except when stated otherwise. Note that in the case when $\mathbb{H}=\R^p$, this means $p$ is a function of $n$. Hence, if a constant $c$ does not depend on $n$, then $c$ does not depend on $p$ either. To  handle the finite and infinite-dimensional settings in a streamlined way, we let $\mathcal{J}=\{1,\dots,\dim(\mathbb{H})\}$ when $\dim(\mathbb{H})<\infty$, and we let $\mathcal{J}=\{1,2,\dots\}$ when $\dim(\mathbb{H})=\infty$. Within this basic setup, our main model assumptions are stated below.

\noindent \begin{assumption}[Data generating model]\label{A:model}
The random vectors $X_1,\dots,X_n\in \mathbb{H}$ are centered and i.i.d., with a compact covariance operator $\Sigma=\E(X_1\otimes X_1)$ whose spectral decomposition is $\sum_{j\in\mathcal{J}} \lambda_j(\Sigma)\phi_j\otimes \phi_j$, with $\|\phi_j\|=1$ for all $j\in\mathcal{J}$.
\begin{enumerate}[(a).]
\item There is a parameter $\nu>0$, fixed with respect to $n$, such that the random variables $\{\zeta_{1j}\}_{j\in\mathcal{J}}$ defined by $\zeta_{1j}=\langle X_1,\phi_j\rangle/\|\langle X_1,\phi_j\rangle\|_{L^2}$ satisfy
$$\sup_{j\in\mathcal{J}}\, \sup_{q\geq 2} \,q^{-\nu} \|\zeta_{1j}\|_{L^q} \ \lesssim \ 1.$$

\item There is a parameter $\beta>1/2$, fixed with respect to $n$, such that the eigenvalues of $\Sigma$ satisfy
\begin{equation*}
1\ \lesssim \ \inf_{j\in\mathcal{J}}\lambda_j(\Sigma) j^{2\beta} \ \leq  \ \displaystyle\sup_{j\in\mathcal{J}}\ts \lambda_j(\Sigma)j^{2\beta} \ \lesssim \ 1.
\end{equation*}
 \item For each $d\in\mathcal{J}$,
let $\J_d:\mathbb{S}^{d\times d}\to\mathbb{S}^{d\times d}$ be the covariance operator of $(\zeta_{11},\dots,\zeta_{1d})\otimes (\zeta_{11},\dots,\zeta_{1d})$, viewed as a random element of $\mathbb{S}^{d\times d}$.
 Then,
\begin{equation*}\label{eqn:eigcond}
1 \ \lesssim \ \inf_{d\in\mathcal{J}}\lambda_{\min}(\J_d) \ \leq \ \sup_{d\in\mathcal{J}}\lambda_{\max}(\J_d) \ \lesssim \ 1.
\end{equation*}
\end{enumerate}
\end{assumption}
\noindent\textbf{Remarks.} The random variables $\{\zeta_{1j}\}_{j\in\mathcal{J}}$ are often referred to as the standardized principal component scores of $X_1$, which appear in the Karhunen-Lo\`eve expansion $X_1=\sum_{j\in\mathcal{J}}\sqrt{\lambda_j(\Sigma)}\zeta_{1j}\phi_j$~\citep{Hsing:2015, Bosq, Horvath}. More generally, this expansion serves as a framework for formulating model assumptions in functional data analysis~\citep[e.g.][]{jacques,pace,jirak2023}. 
With regard to the eigenvalues $\{\lambda_j(\Sigma)\}_{j\in\mathcal{J}}$, decay constraints akin to condition (b) are especially common~\citep[e.g.][]{hallhorowitz,CaiHall}, and the condition $\beta>1/2$ is a natural enhancement of the basic inequality $\lambda_j(\Sigma)\leq \tr(\Sigma) j^{-1}$ that holds for all $j\in\mathcal{J}$. The main novelty of Assumption~\ref{A:model} is condition (c), which allows for a variety of high-dimensional and infinite-dimensional models to be handled in a unified way, as will be explained in Section~\ref{sec:examples}.

As one more preliminary item, we must define the bootstrap version of the sample covariance operator. For this purpose, we will focus on the non-parametric (or empirical) bootstrap~\citep{Efron}, which generates $X_1^*,\dots,X_n^*$ by sampling with replacement from $X_1,\dots,X_n$, and then uses the matrix $\hat\Sigma^* = \frac{1}{n}\sum_{i=1}^n X_i^*\otimes X_i^*$ to produce $\sqrt{n}\|\hat\Sigma^*-\hat\Sigma\|\op$ as a bootstrap sample of the original statistic $\sqrt n\|\hat\Sigma-\Sigma\|\op$.

The following theorem is our main result.

\begin{theorem}\label{thm:main}
 Fix any $\e\in (0,1)$, and suppose that Assumption~\ref{A:model} holds. Then, there is a constant $c>0$ not depending on $n$ such that the event
\begin{equation}\label{eqn:thm:main}
\sup_{t\in\R}\bigg|\P\Big(\sqrt n \|\hat{\Sigma}-\Sigma\|\op\leq t\Big) -\P\Big(\sqrt n\|\hat\Sigma^*-\hat\Sigma\|\op\leq t\,\Big|X\Big)\bigg| \ \leq \ c\, n^{-\frac{\beta-1/2}{2\beta+4+\e}}
\end{equation}
holds with probability at least $1-c/n$. 
\end{theorem}

\noindent\textbf{Remarks.} Although the constant $c$ does not depend on $n$, it should be noted that $c$ may depend on $\e$, $\beta$, $\nu$, as well as constants that are implicit in the expression $\lesssim 1$ in Assumption~\ref{A:model}. The proof is outlined in Section~\ref{sec:outline}, and the detailed arguments are given in subsequent sections.  Regarding the quantity $\beta-1/2$ in the exponent of the bound, the numerical experiments in~\citep{Lopes:Bernoulli} provide strong evidence that this quantity is essential. Those experiments demonstrate that if $\mathbb{H}=\R^p$, and if both $p$ and $n$ are large, then there is a sharp transition in the performance of the bootstrap, from accurate to inaccurate, depending on whether $\beta>1/2$ or $\beta<1/2$. Hence, this suggests that the assumption $\beta>1/2$ is unavoidable for ensuring bootstrap consistency in high dimensions, and that the rate in Theorem~\ref{thm:main} behaves correctly in the extreme case when $\beta\approx 1/2$. Meanwhile, for the other extreme case $\beta\approx \infty$, the rate is near $n^{-1/2}$, which appears ideal considering that a $n^{-1/2}$ rate is generally unimprovable in the univariate Berry-Esseen theorem. In the intermediate regime $1/2\ll  \beta\ll \infty$, we leave the problem of determining the optimal rate of bootstrap approximation for future work. Another related open question is whether the rate can be formulated in a way that does not depend on the decay profile $\lambda_j(\Sigma)\asymp j^{-2\beta}$. That is, it would be of interest to express the rate in terms of spectral parameters that can be associated with $\Sigma$ in greater generality. A prototypical example of such a parameter is the effective rank $\tr(\Sigma)/\|\Sigma\|\op$, which naturally respects the scale-invariance of the Kolmogorov metric, and is known to be a key determinant of the magnitude of $\|\hat\Sigma-\Sigma\|\op$~\citep[][]{Koltchinskii:Bernoulli}. Nevertheless, in the context of bootstrap approximation, it is unclear to what extent this or other parameters may be needed.

Looking towards other possible extensions of Theorem~\ref{thm:main}, it would be of interest to determine the extent to which similar rates of approximation can be developed for resampling methods other than the non-parametric bootstrap.
Some examples include ``weighted'' or ``multiplier'' bootstraps, which are often defined so that $\sqrt n\|\hat\Sigma^*-\hat\Sigma\|\op =\frac{1}{\sqrt n}\|\sum_{i=1}^n \xi_i X_i\otimes X_i\|\op$ for certain types of multiplier random variables $\xi_1,\dots,\xi_n$~\citep{Zhou:2018,Naumov:2019,jirakboot,yu2022testing,ding2023extreme}. This reduces to the non-parametric bootstrap when $(\xi_1,\dots,\xi_n)=(M_1-1,\dots,M_n-1)$, with $(M_1,\dots,M_n)$ being drawn from a Multinomial$(n,\frac{1}{n},\dots,\frac{1}{n})$ distribution, independently of the data. Among other possible choices for the multipliers, some well-known options are to generate $\xi_1,\dots,\xi_n$ as i.i.d.~$N(0,1)$ or Rademacher random variables, independently of the data. Beyond these examples, parametric bootstrap methods can also be considered,
which generate $X_1^*,\dots,X_n^*$ from a parametric distribution that has been fitted to the data, and then use $\sqrt{n}\|\hat\Sigma^*-\hat\Sigma\|\op$ with $\hat\Sigma^* = \frac{1}{n}\sum_{i=1}^n X_i^*\otimes X_i^*$. In recent years, parametric bootstrap methods have shown promise in solving inference problems related to high-dimensional covariance matrices, as these methods are able to directly account for certain model structures~\citep{Lopes:Biometrika,LSS_EJS,yu2023testing}.

\section{Examples of conforming models}\label{sec:examples}
This section presents a variety of models that satisfy Assumption~\ref{A:model}, in both the infinite and high-dimensional settings. 

\subsection{Infinite dimensions}
\noindent\textbf{Gaussian models.} In the Hilbert space of square-integrable functions on the unit interval, $\mathbb{H}=L^2[0,1]$, there are several fundamental Gaussian distributions that satisfy all aspects of Assumption~\ref{A:model}. In fact, conditions~\ref{A:model}(a) and (c) hold for any centered Gaussian distribution on $L^2[0,1]$, regardless of the covariance structure. To see this, first note that if $X_1$ is a centered Gaussian random element, then the sequence $\{\zeta_{1j}\}_{j\in\mathcal{J}}$ consists of i.i.d.~$N(0,1)$ random variables, which implies that the moment condition~\ref{A:model}(a) holds with $\nu=1/2$~\cite[][p.~25]{Vershynin:2018}. Also, this property of $\{\zeta_{1j}\}_{j\in\mathcal{J}}$ implies that
 $\mathsf{J}_d$ is equal to twice the identity operator for each $d\in\mathcal{J}$, and hence condition \ref{A:model}(c) holds with $\lambda_{\min}(\mathsf{J}_d)=\lambda_{\max}(\mathsf{J}_d)=2$. (The calculation for verifying this is a special case of equation~\eqref{eqn:JdMP} in the proof of Proposition~\ref{prop:MP} below.) Lastly, the spectrum decay condition~\ref{A:model}(b) is satisfied by several familiar Gaussian processes on [0,1], such as standard versions of Brownian motion, the Brownian bridge, and the Ornstein-Uhlenbeck process, which all correspond to $\beta=1$. Moreoever, by considering integrated Gaussian processes, it is possible to develop examples corresponding to any integer $\beta>1$. For further details about the eigenvalues associated with these processes, we refer to~\citep[][Ch.~4~pp.~70-71, and Prop.~10]{Ritter}. \\

\noindent\textbf{Non-Gaussian models.} Let us return to the general situation where $\mathbb{H}$ is only assumed to be separable. Here, a model can be specified by defining the observations to be of the form $X_1=\sum_{j\in\mathcal{J}} \sqrt{\lambda_j(\Sigma)}\zeta_{1j}\phi_j$ for a sequence of uncorrelated and standardized random variables $\{\zeta_{1j}\}_{j\in\mathcal{J}}$, while allowing $\{\lambda_j(\Sigma)\}_{j\in\mathcal{J}}$ to be any sequence that satisfies condition~\ref{A:model}(b), and allowing  $\{\phi_j\}_{j\in\mathcal{J}}$ to be any orthonormal basis. For condition \ref{A:model}(a) to hold, it is sufficient for the random variables $\{\zeta_{1j}\}_{j\in\mathcal{J}}$ to be i.i.d.~copies of $\zeta_{11}$, where the latter random variable is centered with unit variance and satisfies $\sup_{q\geq 2}q^{-\nu}\|\zeta_{11}\|_{L^q}\lesssim 1$ for some $\nu>0$ that is fixed with respect to $n$. Notably, this moment condition is weaker than $\zeta_{11}$ having a moment generating function, because if $\zeta_{11}$ does have a moment generating function, then the condition holds with $\nu=1$~\citep[][Prop.~2.7.1]{Vershynin:2018}. 
Lastly, condition \ref{A:model}(c) holds if the random variable $\zeta_{11}$ has one additional property, which is that $\var(\zeta_{11}^2)\gtrsim 1$. This conclusion is reached in two steps, by first considering the formula $\llangle A, \mathsf{J}_d(A)\rrangle= 2\|A\|_F^2+(\var(\zeta_{11}^2)-2)\sum_{j=1}^d A_{jj}^2$, which holds for all $A\in\mathbb{S}^{d\times d}$ and all $d\in\mathcal{J}$ when $\{\zeta_{1j}\}$ are i.i.d.~copies of $\zeta_{11}$~\citep[eqn.~9.8.6]{Bai:Silverstein:2010}. In turn, by maximizing and minimizing this formula for $\llangle A, \mathsf{J}_d(A)\rrangle$ over the set of symmetric matrices with $\|A\|_F=1$,  the eigenvalues of $\mathsf{J}_d$ can be bounded away from zero and infinity via the calculations~\eqref{eqn:JdMP} through \eqref{eqn:MPupper} in the proof of Proposition~\ref{prop:MP} below (while taking $d=p$ and $\Phi=I$ in the notation used there).

\subsection{High dimensions}\label{sec:hd}

For data that reside in a Euclidean space $\mathbb{H}=\R^p$, two of the most popular types of models in high-dimensional statistics and random matrix theory are elliptical models and Mar\v{c}enko-Pastur models. Since both types allow for the spectral decay condition~\ref{A:model}(b), it will only be necessary to show that they are compatible with conditions~\ref{A:model}(a) and~(c). For simplicity, we only define these models in the case of centered random vectors. It is also important to note that our high-dimensional examples do not constrain the size of $p$ with respect to $n$.\\

\noindent\textbf{Elliptical models.} A random vector $X_1\in\R^p$  is said to follow an elliptical model if there exists a non-random positive semidefinite matrix $M\in\mathbb{S}^{p\times p}$ and a random vector $Z_1\in\R^p$ such that
$X_1 \ \overset{}{=}\  MZ_1$,
where $Z_1$ has the properties that  $\E(Z_1Z_1\ttop)=I$, and $QZ_1\overset{{\scriptsize\mathcal{L}}}{=} Z_1$ for any non-random orthogonal matrix $Q\in\R^{p\times p}$. More extensive background on elliptical models can be found in~\citep{Anderson,Muirhead}.

\begin{proposition}\label{prop:elliptical}
Suppose the random vector $X_1\in\R^p$  follows an elliptical model with a positive definite covariance matrix. Then, condition~\ref{A:model}(a) holds when $\sup_{q\geq 2}q^{-\nu}\|Z_{11}\|_{L^q}\lesssim 1$ for some constant $\nu>0$ not depending on $n$, and condition \ref{A:model}(c) holds when $\P(Z_1=0)=0$ and $\var(\|Z_1\|_2^2)\asymp\E(\|Z_1\|_2^2)$.
\end{proposition}

\noindent\textbf{Remarks.} The proof is given in Section~\ref{sec:prop:elliptical}. The condition $\var(\|Z_1\|_2^2)\asymp \E(\|Z_1\|_2^2)$ appears elsewhere in connection with high-dimensional sample covariance matrices, such as in \citep{AOS:elliptical,LSS_EJS}, where it is used to analyze the limit laws of linear spectral statistics in elliptical models. To provide some concrete examples for which the conditions of Proposition~\ref{prop:elliptical} hold, consider a random vector of the form $Z_1=\eta_1 U_1$, where $U_1\in\R^p$ is uniformly distributed on the unit sphere, and $\eta_1$ is a positive random variable that is independent of $U_1$. If the square of $\eta_1$ can be represented as $\eta_1^2=\sum_{j=1}^p \eta_{1j}^2$ for some positive i.i.d.~random variables $\eta_{11},\dots,\eta_{1p}$  satisfying $\E(\eta_{11}^2)=1$, $\var(\eta_{11}^2)\gtrsim 1$, and $\sup_{q\geq 2}q^{-\tilde\nu}\|\eta_{11}^2\|_{L^q}\lesssim 1$ for some fixed $\tilde\nu>0$ that does not depend on $n$, then $Z_1$ satisfies all the conditions in the proposition with $\nu=(\tilde\nu+1)/2$. This class of distributions for $\eta_1^2$ includes familiar parametric distributions such as Chi-Squared$(p)$, Gamma$(p,1)$, and $\tau\cdot\,$Negative-Binomial$(p,\tau)$ for any fixed $\tau\in(0,1)$.

More generally, it is natural to look for examples of conforming elliptical models when $\mathbb{H}$ is infinite-dimensional. (In the context of an abstract separable Hilbert space, a random element is said to be elliptical if all of its one-dimensional projections are elliptical random variables, as defined above for $\R^p$ with $p=1$~\citep[][Definition 2.2]{Tyler:2014}.) However, it turns out that when $\mathbb{H}$ is infinite-dimensional, the class elliptical distributions is much more restricted in a certain sense as compared to the finite-dimensional case. In particular, if $X_1$ is elliptical and has a covariance operator with infinite rank, then $X_1$ must arise from a scale mixture of Gaussian distributions~\citep[][Prop.~2.1]{Tyler:2014}. That is, $X_1$ must be representable in the form  $X_1=\eta_1 V_1$, where $V_1\in\mathbb{H}$ is a centered Gaussian random element, and $\eta_1$ is a non-negative random variable that is independent of $V_1$. This contrasts sharply with the finite-dimensional setting, where many elliptical distributions are not of this type, such as the uniform distribution on the unit sphere of $\R^p$. Also, it is possible to check that if $\mathbb{H}$ is infinite-dimensional, then a scale mixture of Gaussian distributions can only satisfy Assumption~\ref{A:model}(c) when the mixing variable $\eta_1$ is a constant, resulting in a Gaussian distribution.  (See the remark following the proof of Proposition~\ref{sec:prop:elliptical} below for details.) \\

\noindent\textbf{Mar\v{c}enko-Pastur models.} A random vector $X_1\in\R^p$ is said to follow a Mar\v{c}enko-Pastur model if the following holds: There is a non-random positive semidefinite  matrix $M\in\mathbb{S}^{p\times p}$ and a random vector $Z_1\in\R^p$ having i.i.d.~entries with $\E(Z_{11})=0$ and $\var(Z_{11})=1$,  such that
$ X_1\ \overset{}{=} \ M Z_1$.
We refer to~\citep[][]{Bai:Silverstein:2010,Couillet} for additional details on these models, which are sometimes referred to as ``separable'' or ``independent component'' models.

\begin{proposition}\label{prop:MP}
Suppose the random vector $X_1\in\R^p$  follows a Mar\v{c}enko-Pastur  model with a positive definite covariance matrix. Then, condition~\ref{A:model}(a) holds when $\sup_{q\geq 2}q^{-\nu}\|Z_{11}\|_{L^q}\lesssim 1$ for some constant $\nu>0$ not depending on $n$, and condition~\ref{A:model}(c) holds when $\var(Z_{11}^2)\asymp 1$.
\end{proposition}

\noindent\textbf{Remarks.} The proof is given in Section~\ref{sec:prop:MP}. To see that the condition $\var(Z_{11}^2)\asymp 1$ is quite mild, note that the lower bound $\var(Z_{11}^2)\gtrsim 1$ serves to avoid the extreme case when $Z_{11}$ is a Rademacher variable.  That is, if $Z_{11}$ is standardized and $\var(Z_{11}^2)=0$, then $Z_{11}$ must take the values $\pm 1$ with equal probability.

\subsection{Proof of Proposition~\ref{prop:elliptical}}\label{sec:prop:elliptical}
To start, it is helpful to recognize that the matrix $M$ is the same as $\Sigma^{1/2}$. Also, for any $d\in\{1,\dots,p\}$, let $\Phi\in\R^{p\times d}$ be the matrix whose columns are the leading $d$ eigenvectors $\phi_1,\dots,\phi_d$ of $\Sigma$, and let $D=\textsf{Diag}(\lambda_1(\Sigma),\dots,\lambda_d(\Sigma))$. By the definition of $(\zeta_{11},\dots,\zeta_{1d})$ in Assumption~\ref{A:model}, multiplying both sides of the equation $X_1=M Z_1$ with the matrix $D^{-1/2}\Phi\ttop$ implies
\begin{equation}\label{eqn:ellipreln}
(\zeta_{11},\dots,\zeta_{1d}) = \Phi\ttop Z_1.
\end{equation}
Furthermore, due to the orthogonal invariance of $Z_1$, it follows that $\Phi\ttop Z_1 \overset{{\scriptsize\mathcal{L}}}{=} (Z_{11},\dots,Z_{1d})$, and that all the coordinates of $Z_1$ have the same distribution. Therefore, condition~\ref{A:model}(a) is implied by $\sup_{q\geq 2}q^{-\nu}\|Z_{11}\|_{L^q}\lesssim 1$.

 To address condition~\ref{A:model}(c), the relation~\eqref{eqn:ellipreln} implies that $\J_d$ is also the covariance operator of $(\Phi\ttop Z_1)(\Phi\ttop Z_1)\ttop$, and so the formula
\begin{equation}\label{eqn:varellipformula0}
\begin{split}
 \llangle A,\J_d(A)\rrangle   & \ = \  \ \var(\langle \Phi\ttop Z_1, A \Phi\ttop Z_1\rangle)
\end{split}
\end{equation}
holds for any $A\in\mathbb{S}^{d\times d}$. Next, the orthogonal invariance of $Z_1$ with $\P(Z_1=0)=0$ implies that the random variable $\|Z_1\|_2$ and the random vector $Z_1/\|Z_1\|_2$ are independent, with $Z_1/\|Z_1\|_2$ being uniformly distributed on the unit sphere of $\R^p$~\citep[Theorem 1.5.6]{Muirhead}. Based on this, if we let $r_p=\frac{1}{p(p+2)}\E\|Z_1\|_2^4$, and calculate the variance in~\eqref{eqn:varellipformula0} using~\citep[Lemma A.1]{AOS:elliptical}, then it follows that
\begin{equation}\label{eqn:varellipformula}
\begin{split}
  \llangle A,\J_d(A)\rrangle  & \ = \ (r_p-1)\tr(\Phi A\Phi\ttop )^2 \ + \  2r_p\|\Phi A\Phi\ttop \|_F^2,\\
 & \ = \ (r_p-1)\tr(A)^2 \ + \  2r_p\|A\|_F^2.
\end{split}
\end{equation}
We will now establish~2.1\eqref{eqn:eigcond} by using the formula above to separately show that $\lambda_{\max}(\J_d)\lesssim 1$ and $\lambda_{\min}(\J_d)\gtrsim 1$ (with these bounds holding uniformly with respect to $d$). For the largest eigenvalue, we have
\begin{equation*}
\begin{split}
\lambda_{\max}(\J_d) & \ = \ \sup_{\|A\|_F=1} \llangle A, \J_d(A)\rrangle\\[0.2cm]
%
%
& \ \leq \   [(r_p-1)\vee 0]\,d \, + \,  2r_p,
\end{split}
\end{equation*}
where the last step follows from the general inequality $\tr(A)^2 \leq \|A\|_F^2 d$.
 In the case when $r_p<1$, we have $\lambda_{\max}(\J_d)<  2$, and so we may consider the case $r_p\geq 1$, where it is enough to bound $(r_p-1)d+2r_p$.
 For this purpose, observe that the assumptions of Proposition~\ref{prop:elliptical} ensure there is a constant $c>0$ not depending on $p$ such that
\begin{equation}\label{eqn:4thupper}
\E\|Z_1\|_2^4 \ = \  (\E\|Z_1\|_2^2)^2 \ + \ \var(\|Z_1\|_2^2) \ \leq  \ p^2 \ + \ cp,
\end{equation}
where we have used $\E(Z_1Z_1\ttop)=I$. Hence,
\begin{equation*}
\begin{split}
(r_p-1)d +  2r_p & \ = \  \Big(\ts\frac{\E\|Z_1\|_2^4-p(p+2)}{p(p+2)}\Big)d \ + \ \frac{2\E\|Z_1\|_2^4}{p(p+2)}\\[0.2cm]
& \ \leq  \ \Big(\ts\frac{(c-2)p}{p(p+2)}\Big)d \ + \ 2\frac{(p^2+cp)}{p(p+2)}\\[0.2cm]
& \ \lesssim \ 1.
\end{split}
\end{equation*}
Therefore, we conclude that $\lambda_{\max}(\mathsf{J}_d)\lesssim 1$.\\

To lower bound $\lambda_{\min}(\J_d)$, we proceed along similar lines,
\begin{equation*}
\begin{split}
\lambda_{\min}(\J_d)  & \ = \ \inf_{\|A\|_F=1} \llangle A, \J_d(A)\rrangle\\[0.2cm]
& \  = \  \inf_{\|A\|_F=1} (r_p-1)\tr(A)^2+2r_p.
\end{split}
\end{equation*}
In the case when $r_p\geq 1$, this implies $\lambda_{\min}(\J_d)\geq 2$, and so it remains to consider the case when $r_p<1$. This leads to
\begin{equation*}
\begin{split}
\lambda_{\min}(\J_d)  
& \  = \   -\sup_{\|A\|_F=1}\Big((1-r_p)\tr(A)^2\Big) \ + \ 2r_p\\[0.2cm]
& \ = \ -(1-r_p)d \ + \ 2r_p\\[0.2cm]
& \ = \ \Big(\ts\frac{\E\|Z_1\|_2^4 - p(p+2)}{p(p+2)}\Big)d\ + \ 2\ts\frac{\E\|Z_1\|_2^4}{p(p+2)}.
\end{split}
\end{equation*}
Similarly to~\eqref{eqn:4thupper}, we have $\E\|Z_1\|_2^4  \geq p^2+c'p$ for some constant $c'>0$ not depending on $p$. Noting that $d\leq p$, and that the condition $r_p<1$ implies $c'-2<0$, we have
\begin{equation*}
\begin{split}
\lambda_{\min}(\J_d)  &  \ \geq \ \ts\frac{(c'-2)p}{p(p+2)}d \ + \ \ts\frac{2(p^2+c'p)}{p(p+2)}\\[0.2cm]
& \ \geq \ \ts\frac{c'p}{p+2}\\[0.2cm]
& \ \gtrsim \ 1,
\end{split}
\end{equation*}
which completes the proof.\qed

~\\

\noindent\textbf{Remark.} Some steps in the previous proof shed light on Assumption~\ref{A:model}(c) in the case when  $X_1$ is has a centered elliptical distribution and $\textup{dim}(\mathbb{H})=\infty$. Here, we show that if the covariance operator $\Sigma$ of $X_1$ has infinite rank in this context, then $X_1$ can only satisfy Assumption~\ref{A:model}(c) when it is Gaussian. To begin, recall from our discussion in Section~\ref{sec:hd} that the infinite rank of $\Sigma$ requires $X_1$ to be a scale mixture of Gaussian distributions~\citep[][Prop.~2.1]{Tyler:2014}. That is,  $X_1=\eta_1 V_1$, where $V_1\in\mathbb{H}$ is a centered Gaussian random element with covariance operator $\Sigma$, and $\eta_1$ is a non-negative random variable that is independent of $V_1$ with $\E(\eta_1^2)=1$. From this representation of $X_1$, it follows that for any integer $d\geq 1$, we have $(\zeta_{11},\dots,\zeta_{1d})=\eta_1 W_1$, where $W_1\in\R^d$ is a standard Gaussian vector that is independent of $\eta_1$.
Then, for any symmetric matrix $A\in\R^{d\times d}$, a direct calculation gives
\begin{equation*}
\begin{split}
\llangle A,\mathsf{J}_d(A)\rrangle & \ = \ \textup{var}(\eta_1^2 \,\langle  W_1,A W_1\rangle)\\
&  \ = \  \tr(A)^2\var(\eta_1^2)+2\|A\|_F^2(\var(\eta_1^2)+1),
\end{split}
\end{equation*}
and so
\begin{equation*}
\begin{split}
\lambda_{\max}(\mathsf{J}_d) &\ = \  \sup_{\|A\|_F=1}\llangle A,\mathsf{J}_d(A)\rrangle \  = \  (d+2)\var(\eta_1^2)+2.
%
%
\end{split}
\end{equation*}
Thus, the only way $\sup_{d\in\mathcal{J}}\lambda_{\max}(\J_d)$ can be finite in Assumption~\ref{A:model}(c)  is if $\var(\eta_1^2)=0$, which means that $\eta=1$ almost surely, and so $X_1$ must be Gaussian.

\subsection{Proof of Proposition~\ref{prop:MP}} \label{sec:prop:MP}
As in the previous proof, note that the matrix $M$ is equal to $\Sigma^{1/2}$, and for any $d\in\{1,\dots,p\}$, let  $\Phi\in\R^{p\times d}$ be the matrix whose columns are the leading $d$ eigenvectors $\phi_1,\dots,\phi_d$ of $\Sigma$. It follows that
 \begin{equation}\label{eqn:secondzz}
 (\zeta_{11},\dots,\zeta_{1d}) = \Phi\ttop Z_1 ,
  \end{equation}
with the entries on the left side being as defined in Assumption~\ref{A:model}. Due to Rosenthal's inequality (Lemma~\ref{lem:rosenthal} in Section~\ref{sec:background}), we have the following bound for any $q\geq 2$ and $j=1,\dots,d$,
\begin{equation*}
\begin{split}
\|\zeta_{1j}\|_{L^q} & \ = \  \big\| \ts\sum_{l=1}^p \phi_{jl}Z_{1l}\big\|_{L^q}\\
 &\ \lesssim  \  q\cdot\max\bigg\{ \big\|\ts\sum_{l=1}^p \phi_{jl}Z_{1l}\big\|_{L^2} \, , \, \big(\ts\sum_{l=1}^p \big\|\phi_{jl}Z_{1l}\|_{L^q}^q\big)^{1/q}\bigg\}\\
 & \lesssim \  q\cdot\max\big\{1\, , \ q^{\nu}\big(\ts\sum_{l=1}^p |\phi_{jl}|^q \big)^{1/q}\Big\}\\
 &\lesssim \ q^{\nu+1},
\end{split}
\end{equation*}
where we have used the fact that $\phi_j$ is a unit vector. This handles the condition~\ref{A:model}(a), since the exponent $\nu+1$ here does not need to match the exponent $\nu$ appearing in condition~\ref{A:model}(a).

Turning our attention to  condition~\ref{A:model}(c), fix any $A\in\mathbb{S}^{d\times d}$ and observe that~\eqref{eqn:secondzz} gives
\begin{equation*}
\begin{split}
 \llangle A,\J_d(A)\rrangle  
  & \ = \ \var(\langle \Phi\ttop Z_1, A \,\Phi\ttop Z_1\rangle).
\end{split}
\end{equation*}
Letting $\kappa=\E(Z_{11}^4)$ and noting that $\Phi$ has orthogonal columns, we may evaluate the last expression using a known formula for the variance of a quadratic form involving centered i.i.d.~random variables with unit variance~\citep[eqn.~9.8.6]{Bai:Silverstein:2010},
\begin{equation}\label{eqn:JdMP}
\begin{split}
 \llangle A,\J_d(A)\rrangle  
  & \ = \ 2\|A\|_F^2 + \ (\kappa-3)\sum_{j=1}^p [\Phi A \Phi\ttop]_{jj}^2.
  \end{split}
  \end{equation}
Now, we can compute the smallest eigenvalue of $\J_d$ as
\begin{equation*}
\begin{split}
\lambda_{\min}(\mathsf{J}_d) & \ = \ \inf_{\|A\|_F=1} \llangle A\, ,\,\J_d(A)\rrangle\\[0.2cm]
& \ = \ 2 \ + \ \inf_{\|A\|_F=1} (\kappa-3)\Big\llangle \Phi A\Phi\ttop ,\textsf{Diag}(\Phi A\Phi\ttop)\Big\rrangle.
\end{split}
\end{equation*}
In the case when $\kappa\geq 3$, the last expression gives $\lambda_{\min}(\mathsf{J}_d)\geq 2$. To handle the case when $\kappa<3$, note that if $A,B\in\mathbb{S}^{d\times d}$, then $\llangle B, \textsf{Diag}(B)\rrangle \leq \|B\|_F^2$, and $\|\Phi A\Phi\ttop\|_F^2= \|A\|_F^2$. Therefore,
\begin{equation}\label{eqn:MPlower}
\begin{split}
\lambda_{\min}(\mathsf{J}_d) &  \ \geq \ 2 \, + \, (\kappa-3)
 \ = \ \var(Z_{11}^2)
 \ \gtrsim \ 1.
\end{split}
\end{equation}
Similarly, for the largest eigenvalue, we have
\begin{equation}\label{eqn:MPupper}
\begin{split}
\lambda_{\max}(\mathsf{J}_d) & \ = \ \sup_{\|A\|_F=1} \llangle A, \J_d(A)\rrangle\\[0.2cm]
& \ = \ 2 \ + \ \sup_{\|A\|_F=1} (\kappa-3)\Big\llangle \Phi A\Phi\ttop\, ,  \textsf{Diag}(\Phi A \Phi\ttop)\Big\rrangle\\[0.2cm]
& \ \leq \ 2 \, + \, |\kappa-3|\\[0.2cm]
& \ \lesssim \ 1,
\end{split}
\end{equation}
as needed.\qed

\section{Outline of the proof of Theorem~\ref{thm:main}}\label{sec:outline}
Before organizing the main parts of the proof, we need to define a collection of objects that will appear frequently in the arguments.\\

\noindent\textbf{Special indices.} First, let $\e$ be as chosen in the statement of Theorem~\ref{thm:main},
 and define the associated integer
\begin{equation}\label{eqn:mdef}
m=\Big\lceil n^{\frac{\e(\beta-1/2)}{12\beta(2\beta+5)^2}}\wedge \dim(\mathbb{H})\Big\rceil.
\end{equation}
The exponent on $n$ in this definition serves only as a technical expedient, and a smaller choice of the exponent could also be used.
 Second, let $C>2\beta/(2\beta-1)$ be a constant fixed with respect to $n$, and define
\begin{align}
\ell & \ = \ \big\lceil m^{\frac{1}{C}}\big\rceil\label{eqn:elldef}.
\end{align}
Third, define
\begin{align}
k& \ = \  \Big\lceil n^{\frac{1}{2\beta+4}}\wedge \dim(\mathbb{H})\Big\rceil,\label{eqn:kdef}
\end{align}
and note that $m$, $\ell$, and $k$ satisfy the ordering relation $1\leq \ell \leq m\leq k\leq \dim(\mathbb{H})$. \\

\noindent\textbf{Projected data and covariance matrices.} Define the random vectors $Y_1,\dots,Y_n\in\R^k$ according to
\begin{equation*}
Y_i= (\langle X_i,\phi_1\rangle,\dots,\langle X_i,\phi_k\rangle).
\end{equation*}
The common covariance matrix of $Y_1,\dots,Y_n$ is denoted as $\Sigma_k=\E(Y_1Y_1\ttop)\in\mathbb{S}^{k\times k}$, and can alternatively be written as $\Sigma_k=\textsf{Diag}(\lambda_1(\Sigma),\dots,\lambda_k(\Sigma))$.  To lighten notation, we will also use $\zeta_1\in\R^k$ to denote the random vector 
\begin{equation}\label{eqn:zetadef}
\zeta_1=(\zeta_{11},\dots,\zeta_{1k}),
\end{equation}
with entries as defined in Assumption~\ref{A:model}.

The empirical and bootstrap versions of $\Sigma_k$ are respectively denoted as
\begin{align}
\hat\Sigma_k & \ = \ \ts\frac 1n \sum_{i=1}^n Y_i Y_i\ttop\\[0.2cm]
 \hat\Sigma_k^* & \ = \ \ts\frac{1}{n}\sum_{i=1}^n Y_i^*(Y_i^*)\ttop,
\end{align}
where $Y_1^*,\dots,Y_n^*$ are sampled with replacement from $Y_1,\dots,Y_n$.\label{Ystardef}\\

\noindent\textbf{Gaussian random matrices.} 
Let $G$ be a centered Gaussian random matrix in $\mathbb{S}^{k\times k}$ with the same covariance operator as $\sqrt n (\hat\Sigma_k-\Sigma_k)$. (This also means that $G$ has the same covariance operator as $Y_1Y_1\ttop$.) The bootstrap version of $G$ is denoted as $G^*$, which is defined to be a random matrix in $\mathbb{S}^{k\times k}$ whose conditional distribution $\mathcal{L}(G^*|X)$ is centered and Gaussian with the same covariance operator as $\mathcal{L}(\sqrt n(\hat\Sigma_k^*-\hat\Sigma_k)|X)$.\\

\noindent\textbf{Subsets of the ball.}
Let $\mathbb{B}^k$ denote the unit ball for the Euclidean norm in $\R^k$, and define the subsets
\begin{align}
\mathbb{B}_{\scaleto{\!\triangle}{4pt}}^k & \ = \ \Big\{u\in\mathbb{B}^k \ \Big| \ \ts\sum_{j=1}^m \langle u,e_j\rangle^2 \, \geq \, m^{-2\beta+1}\Big\}\label{eqn:Bktriangledef}\\[0.3cm]
 \mathbb{B}_{\scaleto{\!\triangledown}{4pt}}^k\! & \ = \ \mathbb{B}^k \setminus \mathbb{B}_{\scaleto{\!\triangle}{4pt}}^k.
\end{align}
For any $A\in\mathbb{S}^{k\times k}$, define the associated quantities
\begin{align}
 \|A\|_{_{\scaleto{\!\triangle}{4pt}}} & \ = \ \sup_{u\in \mathbb{B}_{\scaleto{\!\triangle}{4pt}}^k} |\langle u, A u\rangle|\label{eqn:triupdef}\\[0.3cm]
 \|A\|_{\scaleto{\!\triangledown}{4pt}}& \ = \ \sup_{u\in \mathbb{B}_{\scaleto{\!\triangledown}{4pt}}^k} |\langle u, A u\rangle|.
\end{align}
The intuition underlying the definition~\eqref{eqn:Bktriangledef} is that  $\mathbb{B}_{\scaleto{\!\triangle}{4pt}}^k$ consists of vectors in $\mathbb{B}^k$ that are ``suitably aligned'' with the first $m$ standard basis vectors. For the  random matrix $G$, these directions play a dominant role, in the sense that the random variable $\|G\|{_{\scaleto{\!\triangle}{4pt}}}$ serves as an approximation to $\|G\|\op$. Moreover, because $\mathbb{B}_{\scaleto{\!\triangle}{4pt}}^k$ excludes vectors $u$ for which the variance of $\langle u, G u\rangle$ is very small, it becomes feasible to establish anti-concentration properties of $\|G\|{_{\scaleto{\!\triangle}{4pt}}}$, which can then be transferred to $\|G\|\op$ and other random variables in our analysis. Lastly, the reason for introducing $ \|\cdot\|_{\scaleto{\!\triangledown}{4pt}}$ is that the problem of showing that $\|G\|{_{\scaleto{\!\triangle}{4pt}}}$ and $\|G\|\op$ are closely coupled can be partially reduced to the problem of showing that $ \|G\|_{\scaleto{\!\triangledown}{4pt}}$ is likely to be small.\\

\noindent\textbf{A decomposition for Theorem~\ref{thm:main}.} Recall that we use $d_{\textup{K}}$ to denote the Kolmogorov metric, as described in the notation paragraph of Section~\ref{sec:intro}. The proof is organized by bounding the distance of interest with five terms,
$$d_{\textup{K}} \bigg(\mathcal{L}\Big(\sqrt n \|\hat{\Sigma}-\Sigma\|\op\Big) \, , \, \mathcal{L}\Big(\sqrt n\|\hat\Sigma^*-\hat\Sigma\|\op\Big|X\Big)\bigg) \ \ \leq \ \ \I_n \   + \ \II_n  \ + \ \widehat{\III}_n \ + \ \hat{\II}_n  \ + \ \hat{\I}_n.$$
The definitions of the five terms, with brief descriptors, are as follows. 
\begin{align*}\label{romandefs}
\shortintertext{Dimension reduction:}\\[-1.25cm]
 \ \ \ \ \ \ \ \ \ \ \ \  \ \ \ \ \ \ \ \ \ \  \ \ \ \  \ \ \  \ \ \  \ \ \ \ \ \ \ \ \ \ \ \ \ \   \ \ \ \ \ \I_n & \ = \ d_{\textup{K}}\Big(\mathcal{L}\big(\sqrt n \|\hat{\Sigma}-\Sigma\|\op\big) \, , \, \mathcal{L}\big(\sqrt n \|\hat\Sigma_k-\Sigma_k\|\op\big)\Big)\\[1cm]
\shortintertext{Gaussian approximation:}\\[-1.25cm]
\II_n& \ = \ d_{\textup{K}}\Big(\mathcal{L}\big( \sqrt n \|\hat\Sigma_k-\Sigma_k\|\op\big) \, , \, \mathcal{L}\big( \|G\|\op \big)\Big)\\[1cm]
\shortintertext{Gaussian comparison:}\\[-1.25cm]
\widehat{\III}_n& \ = \ d_{\textup{K}}\Big(\mathcal{L}\big(\|G\|\op \big),\mathcal{L}\big( \|G^*\|\op \big|X\big)\Big)\\[1cm]
\shortintertext{Bootstrap approximation:}\\[-1.25cm]
\hat{\II}_n& \ = \ d_{\textup{K}}\Big(\mathcal{L}\big(\|G^*\|\op\big|X\big) \ , \ \mathcal{L}\big(\sqrt n\|\hat\Sigma_k^*-\hat\Sigma_k\|\op \big|X\big)\Big)\\[1cm]
\shortintertext{Dimension reduction:}\\[-1.25cm]
\hat{\I}_n& \ = \ d_{\textup{K}}\Big(\mathcal{L}\big(\sqrt n\|\hat\Sigma_k^*-\hat\Sigma_k\|\op\big|X\big) \  , \ \mathcal{L}\big(\sqrt n\|\hat\Sigma^*-\hat\Sigma\|\op\big|X\big)\Big)
\end{align*}
\normalsize
Bounds on the five terms are developed in the respective Sections~\ref{sec:I},~\ref{sec:III},~\ref{sec:hatIV},~\ref{sec:III}, and~\ref{sec:hatI}. Once the bounds resulting from these sections are combined, the proof of Theorem~\ref{thm:main} is complete.\qed

\section{The term $\I_n$: Dimension reduction}\label{sec:I}
\begin{proposition}\label{prop:I}
Suppose the conditions of Theorem~\ref{thm:main} hold. Then,
\begin{equation*}
\I_n \ \lesssim \ n^{-\frac{\beta-1/2}{2\beta+4+\e}}.
\end{equation*}
\end{proposition}

\proof The Kolmogorov distance between the distributions of any two random variables $U$ and $V$ on the same probability space can be bounded as
\begin{equation}\label{eqn:coupleanti}
\dK(\LL(V),\LL(U)) \ \leq \  \P(|V-U|\geq \delta)  \ + \ \sup_{t\in\R}\P\big(|U-t|\leq \delta\big),
\end{equation}
for any $\delta>0$.
For the remainder of the current proof, we will take $V=\sqrt n \|\hat{\Sigma}-\Sigma\|\op$ and $U=\sqrt n\|\hat\Sigma_k-\Sigma_k\|\op$ so that $\I_n=\dK(\LL(V),\LL(U))$. Also note that we may assume without loss of generality that $k<\dim(\mathbb{H})$, for otherwise $V=U$, and then $\textup{I}_n=0$.

The coupling probability $ \P(|V-U|\geq \delta)$ is handled by Lemma~\ref{lem:covcouple} below, which ensures that there is a constant $c>0$ not depending on $n$ such that the following event holds with probability at most $c/n$,
\begin{equation*}
 \bigg|\sqrt n \|\hat{\Sigma}-\Sigma\|\op- \sqrt n\|\hat\Sigma_k-\Sigma_k\|\op\bigg| \ \geq \  c\,  n^{-\frac{\beta-1/2}{2\beta+4}}\log(n)^{\nu+1/2}.
\end{equation*}
Accordingly, we will take $\delta= c\,  n^{-\frac{\beta-1/2}{2\beta+4}}\log(n)^{\nu+1/2}$ in~\eqref{eqn:coupleanti}.

To handle the anti-concentration probability in~\eqref{eqn:coupleanti}, we will use the fact that if two random variables are close in the Kolmogorov metric, then they have similar anti-concentration behavior. More precisely, if $W$ is any other random variable, then
\begin{equation}\label{eqn:reversecoupleanti}
 \sup_{t\in\R}\P\big(|U-t|\leq \delta \big)  \ \leq \ \sup_{t\in\R}\P\big(|W-t|\leq \delta\big)  \ + \ 2\dK(\LL(U),\LL(W)).
\end{equation}
In the current context, we will take $W=\|G\|_{_{\scaleto{\!\triangle}{4pt}}}$,  and Lemma~\ref{lem:Gaussiananti} will show that $\sup_{t\in\R}\P\big(|\|G\|_{_{\scaleto{\!\triangle}{4pt}}}-t|\leq \delta\big)$ is at most of order $n^{-\frac{\beta-1/2}{2\beta+4+\e}}$ based on the choice of $\delta$ mentioned earlier. Hence, it remains to bound $\dK(\LL(U),\LL(W))$, which will be based on a comparison with $\|G\|\op$, namely
\begin{equation*}
\dK(\LL(U),\LL(W)) \ \leq \ \dK(\LL(U), \LL(\|G\|\op)) \ + \ \dK(\LL(\|G\|\op), \LL(W)).
\end{equation*}
The first term on the right is the same as $\II_n$, and Proposition~\ref{prop:optriI} will show that it is at most of order $n^{-\frac{\beta-1/2}{2\beta+4+\e}}$. Lastly, Proposition~\ref{prop:optri} will show that $\dK(\LL(\|G\|\op), \LL(W))$ is at most of order $1/n$, which completes the proof.\qed

\begin{lemma}\label{lem:covcouple}
Suppose the conditions of Theorem~\ref{thm:main} hold. Then, there is a constant $c>0$ not depending on $n$ such that the following event holds with probability at most $c/n$,
\begin{equation*}
 \bigg|\sqrt n \|\hat{\Sigma}-\Sigma\|\op- \sqrt n\|\hat\Sigma_k-\Sigma_k\|\op\bigg| \ \geq \  c  n^{-\frac{\beta-1/2}{2\beta+4}}\log(n)^{\nu+1/2}.
\end{equation*}
\end{lemma}

\noindent \textbf{Remark.} Strictly speaking, the difference $\sqrt n\|\hat\Sigma-\Sigma\|\op-\sqrt n\|\hat\Sigma_k-\Sigma_k\|\op$ is always non-negative, but we have retained the absolute value, since it makes the statement of the lemma appear more natural when it is applied in the proof of Proposition~\ref{prop:I}.

\proof Define the projection operator $\mathsf{\Pi}_k:\mathbb{H}\to\mathbb{H}$ by
$$\mathsf{\Pi}_k(x) = \sum_{j=1}^k\langle x,\phi_j\rangle \phi_j,$$
and let $\mathsf I$ denote the identity operator on $\mathbb{H}$.
Also, recall the notations $\hat\Sigma_k=\frac 1n \sum_{i=1}^nY_iY_i\ttop$ and $\Sigma_k=\E(\hat\Sigma_k)$, where $Y_i=(\langle X_i,\phi_1\rangle,\dots,\langle X_i,\phi_k\rangle)$ for all $i=1,\dots,n$.  Since the map $x\mapsto (\langle x,\phi_1\rangle,\dots,\langle x,\phi_k\rangle)$ is an isometry from $\textup{span}\{\phi_1,\dots,\phi_k\}$ to $\R^k$, the variational representation for the operator norm of a compact operator~\cite[Theorem 4.3.3]{Hsing:2015} implies
$$\|\hat\Sigma_k-\Sigma_k\|\op \ = \  \|\mathsf\Pi_k(\hat\Sigma-\Sigma)\mathsf\Pi_k\|\op.$$ 
Therefore,
\begin{equation}\label{eqn:splitnorms}
\begin{split}
\bigg|\|\hat{\Sigma}-\Sigma\|\op- \|\hat\Sigma_k-\Sigma_k\|\op\bigg|
& \ \leq \ \Big\|(\hat\Sigma-\Sigma)\, -\, \mathsf\Pi_k(\hat\Sigma-\Sigma)\mathsf\Pi_k\Big\|\op\\
& \ = \ \Big\| (\mathsf I -\mathsf\Pi_k)(\hat\Sigma -\Sigma) + \mathsf \Pi_k(\hat\Sigma-\Sigma)(\mathsf I-\mathsf\Pi_k)\Big\|\op\\
& \ \leq \ 2\Big\|(\mathsf I -\mathsf\Pi_k)(\hat\Sigma -\Sigma)\Big\|\op\\
& \ = \  2\Bigg\|\frac{1}{n}\sum_{i=1}^n  [(\mathsf I-\mathsf\Pi_k)X_i] \otimes X_i \ - \ \E\big([(\mathsf I-\mathsf\Pi_k)X_i] \otimes X_i \big)\Bigg\|\op.
\end{split}
\end{equation}
\normalsize
Next, let $\mathbb{B}$ be the unit ball of $\mathbb{H}$, and define $\mathbb{B}_{-k}\subset \mathbb{B}$ as the subset of vectors $v\in\mathbb{B}$ satisfying $\langle v,\phi_j\rangle=0$ for all $j=1,\dots,k$.  Hence, for any $u\in\mathbb{B}$, there exists some $v\in\mathbb{B}_{-k}$ such that $\langle u,(\mathsf I-\mathsf\Pi_k)X_i\rangle = \langle v,X_i\rangle$. So, applying the variational representation of the operator norm to the last expression in~\eqref{eqn:splitnorms},  we have
\begin{equation}\label{eqn:prev}
\begin{split}
\bigg|\|\hat{\Sigma}-\Sigma\|\op- \|\hat\Sigma_k-\Sigma_k\|\op\bigg|
 & \ \ \leq \ \displaystyle   \ \sup_{(v,w)\in\mathbb{B}_{-k}\times \mathbb{B}}2 \Bigg|\frac{1}{n}\sum_{i=1}^n \langle v, X_i\rangle \langle X_i,w\rangle \, - \, \E(\langle v, X_i\rangle \langle X_i,w\rangle)\Bigg |.
\end{split}
\end{equation}
\normalsize
Let $t>0$ be a free parameter to be chosen later, and define the vectors
$$\omega(t)=\ts\frac{tv}{2}+\ts\frac{w}{2t} \text{ \ \ \ \ \ and \ \  \ \ \ } \tilde\omega(t) =\frac{tv}{2}-\frac{w}{2t}. $$
Because the bilinear form $(v,w)\mapsto \langle v, X_i\rangle \langle X_i,w\rangle$ takes the same value when acting on $(v,w)$ and $(tv, \frac{w}{t})$, the polarization identity implies 
$$\langle v, X_i\rangle \langle X_i,w\rangle  \ = \ \langle \omega(t),X_i\rangle^2 \ - \  \langle \tilde \omega(t),X_i\rangle^2.$$
To control the squares on the right side, it is helpful to define the linear operator ${\mathsf{A}}(t):\mathbb{H}\to\mathbb{H}\oplus\mathbb{H}$ according to
$${\mathsf{A}}(t)x \ = \  \Big(\ts t(\mathsf{I}-\mathsf{\Pi}_k) x \, , \,  \frac{x}{t} \Big).$$
Based on this definition, the following relations hold simultaneously for all $t>0$ when the vectors $v$ and $w$ satisfy $(v,w)\in \mathbb{B}_{-k}\times\mathbb{B}$, 
$$\langle \omega(t),X_i\rangle \ = \ \Big\langle (\ts\frac{v}{2},\ts\frac{w}{2}) \, , \, {\mathsf{A}}(t)X_i\Big\rangle \text{ \ \ \ and \ \ \ } \langle \tilde\omega(t),X_i\rangle \ = \ \Big\langle (\ts\frac{v}{2},-\ts\frac{w}{2}) \, , \, {\mathsf{A}}(t)X_i\Big\rangle.$$
Furthermore, by letting  $\xi_i(t)={\mathsf{A}}(t)X_i$, and noting that $(\ts\frac{v}{2},\ts\frac{w}{2})$ and $(\ts\frac{v}{2},-\ts\frac{w}{2})$ lie in the unit ball $\bar{\mathbb{B}}$ of $\mathbb{H}\oplus \mathbb{H}$, it follows from~\eqref{eqn:prev} that
\begin{equation}\label{eqn:Vop}
\begin{split}
\bigg|\|\hat{\Sigma}-\Sigma\|\op- \|\hat\Sigma_k-\Sigma_k\|\op\bigg|
 & \ \leq \ \displaystyle   \ \sup_{u\in \bar{\mathbb{B}}} \,  \,  4\Bigg|\frac{1}{n}\sum_{i=1}^n \langle u,\xi_i(t)\rangle^2 - \, \E \langle u,\xi_i(t)\rangle^2\Bigg |\\[0.2cm]
 &  \ = \  4\Bigg\| \frac{1}{n}\sum_{i=1}^n \xi_i(t)\otimes \xi_i(t) - \E \,\xi_i(t)\otimes  \xi_i(t)\Bigg\|\op.
\end{split}
\end{equation}
\normalsize
The norm in the last line can be handled using Lemma~\ref{lem:op}, which provides a general bound on the operator norm error of sample covariance matrices in separable Hilbert spaces, developed previously in~\citep{Lopes:Bernoulli}. Specifically, if we let $q=\log(n)\vee 3$, then this result gives
\small
\begin{equation}\label{eqn:bigLqbound}
\begin{split}
\Bigg(\E\bigg|\|\hat{\Sigma}-\Sigma\|\op  & - \|\hat\Sigma_k-\Sigma_k\|\op\bigg|^q\Bigg)^{1/q}
%
%
 \ \lesssim \ \  \bigg( \ts\frac{\sqrt q}{\sqrt n}\cdot \|\E \xi_1(t)\otimes \xi_1(t)\|\op^{1/2}\cdot (\E\|\xi_1(t)\|^{2q})^{\frac{1}{2q}}\bigg) \bigvee  \bigg(\ts\frac{q}{n}\cdot(\E \|\xi_1(t)\|^{2q})^{\frac{1}{q}}\bigg).
\end{split}
\end{equation}
\normalsize
To bound the right side, first observe that
\begin{equation}\label{eqn:vecLqbound}
\begin{split}
\|\E \xi_1(t)\otimes \xi_1(t)\|\op^{1/2} 
& \ = \  \|{\mathsf{A}}(t)\Sigma^{1/2}\|\op\\[0.2cm]
& \ \leq \  t\lambda_k(\Sigma)^{1/2} \ + \ \ts\frac{1}{t}\lambda_1(\Sigma)^{1/2}\\[0.2cm]
& \ \lesssim \ tk^{-\beta}  \ + \ \ts\frac{1}{t}.
\end{split}
\end{equation}
Next,  we have
\begin{equation}\label{eqn:vecLqbound2}
\begin{split}
(\E\|\xi_1(t)\|^{2q})^{\frac{1}{2q}} 
& \ \leq  \  t \Big\|  \ts\sum_{j\geq k+1} \langle X_1,\phi_j\rangle^2\Big\|_{L^{q}}^{1/2} \ + \ \ts\frac{1}{t} \Big\| \ts\sum_{j\geq 1} \langle X_1,\phi_j\rangle^2\Big\|_{L^{q}}^{1/2} \\
& \ \leq  \  t\Big(\ts\sum_{j\geq k+1} \|\langle X_1, \phi_j\rangle\|_{L^{2q}}^2\Big)^{1/2} \ +\ \frac{1}{t}\Big(\sum_{j\geq 1}\|\langle X_1, \phi_j\rangle\|_{L^{2q}}^2\Big)^{1/2}\\[0.2cm]
& \ \lesssim \ q^{\nu} t \Big(\ts\sum_{j\geq k+1} j^{-2\beta}\Big)^{1/2} \ + \ \ts\frac{q^{\nu}}{t}\Big(\sum_{j\geq 1} j^{-2\beta}\Big)^{1/2} \ \ \ \ \ \ \text{(Assumptions~\ref{A:model}(a) and (b))}\\[0.2cm]
& \ \lesssim \ q^{\nu}\Big( t k^{-\beta+1/2} \ + \ \ts\frac{1}{t}\Big).
%
\end{split}
\end{equation}
\normalsize
Taking $t=k^{\beta/2}$ in~\eqref{eqn:vecLqbound} and~\eqref{eqn:vecLqbound2} gives
\begin{equation}\label{eqn:EVotimesV}
 \|\E \,\xi_1(t)\otimes \xi_1(t)\|\op^{1/2} \ \lesssim \ k^{-\beta/2},
 \end{equation}
and 
\begin{equation}\label{eqn:V1norm}
(\E\|\xi_1(t)\|^{2q})^{\frac{1}{2q}} \ \lesssim \ q^{\nu}\, k^{-\beta/2+1/2}.
\end{equation}
Combining these bounds with ~\eqref{eqn:bigLqbound}, we have
\begin{equation}\label{eqn:lastdimred}
\begin{split}
 \Bigg(\E\bigg|\|\hat{\Sigma}-\Sigma\|\op  - \|\hat\Sigma_k-\Sigma_k\|\op\bigg|^q\Bigg)^{1/q} 
 & \ \lesssim \ \Big(\ts\frac{q^{\nu+1/2}}{\sqrt n} k^{-\beta+1/2}\Big)\vee \Big(\ts\frac{q^{2\nu+1}}{n}  k^{-\beta+1}\Big)\\[0.2cm]
 & \ \lesssim \  \ts\frac{q^{\nu+1/2} \,k^{-\beta+1/2}}{\sqrt n}.
 \end{split}
\end{equation}
To finish the proof, note that we may assume that $k<\dim(\mathbb{H})$, for otherwise $ \|\hat{\Sigma}-\Sigma\|\op$ and $\|\hat\Sigma_k-\Sigma_k\|\op$ are equal. Consequently, the definition of $k$ in~\eqref{eqn:kdef} implies $k\geq n^{\frac{1}{2\beta+4}}$. Finally, rescaling both sides of~\eqref{eqn:lastdimred} by $\sqrt n$ and applying Chebyshev's inequality leads to the stated result.\qed

\subsection{Anti-concentration inequality for $\|G\|_{_{\scaleto{\!\triangle}{4pt}}}$}

\begin{lemma}\label{lem:Gaussiananti}
Suppose the conditions of Theorem~\ref{thm:main} hold, and let $\delta_n\in (0,1)$ be any numerical sequence that satisfies $\sqrt{\log(1/\delta_n)}\lesssim m^{4\beta-1}$. Then, 
\begin{equation*}
\sup_{t\in\R}\ \P\Big(\big|\|G\|_{_{\scaleto{\!\triangle}{4pt}}}-t\big| \ \leq  \ \delta_n\Big) \ \lesssim \ \delta_n \, m^{12\beta -3}.
\end{equation*}
\end{lemma}
\proof  Define the set $\Theta=\mathbb{B}_{\scaleto{\!\triangle}{4pt}}^k\times \{\pm 1\}$, whose generic element is written as $\theta=(v,s)$. Consider the centered Gaussian process $\Gamma(\theta)= s v\ttop G v$ indexed by $\Theta$, which satisfies
$$\|G\|_{_{\scaleto{\!\triangle}{4pt}}} \ = \ \sup_{\theta\in\Theta}\,\Gamma(\theta).$$
For each $\theta\in\Theta$, define $\varsigma^2(\theta) =\var(\Gamma(\theta))$,
as well as
\begin{equation}\label{eqn:varsigupperlower}
 \bar\varsigma =\sup_{\theta\in\Theta}\varsigma_n(\theta) \ \ \ \ \text{ and } \ \ \ 
\vunderline{\varsigma} =\inf_{\theta\in\Theta} \varsigma_n(\theta).
\end{equation}
In addition, define the expected supremum
$$\mu=\E\Big(\ts\sup_{\theta\in\Theta}\Gamma(\theta)/\varsigma(\theta)\Big).$$
The anti-concentration inequality for Gaussian processes in Theorem 3 of~\citep{CCK:PTRF} implies
\begin{equation}\label{eqn:anticintermed}
\small
\sup_{t\in\R}\P\big(|\|G\|_{_{\scaleto{\!\triangle}{4pt}}}-t|\leq \delta_n\big)  \ \, \lesssim \,  \ 
\ts\frac{\bar{\varsigma}}{ \ \vunderline{\varsigma}^2}\cdot \delta_n\cdot\Big(\mu+\sqrt{1\vee\log(\vunderline{\varsigma}/\delta_n)}\Big).
\end{equation}
(The dependence of this bound on $\bar{\varsigma}$ and $\vunderline{\varsigma}$ is not given explicitly in the quoted result, but a brief inspection of its proof shows that it is sufficient to use a prefactor of $\ts\frac{\bar{\varsigma}}{\vunderline{\varsigma}^2}$.) Lemma~\ref{lem:minvar} below ensures that $\underline \varsigma$ and $\bar\varsigma$ can be controlled as
\begin{equation*}
\underline{\varsigma} \ \gtrsim \ m^{-4\beta+1} \ \ \ \ \text{ and } \ \ \ \ \bar{\varsigma} \ \lesssim \ 1,
\end{equation*}
and so
\begin{equation*}
\ts\frac{\bar{\varsigma}}{\vunderline{\varsigma}^2} \ \lesssim \ m^{8\beta-2} \ \ \ \text{ and }\ \ \  \sqrt{1\vee\log(\vunderline{\varsigma}/\delta_n)} \ \lesssim \ m^{4\beta-1}.
\end{equation*}

Next, we bound $\mu$ as
\begin{equation*}
\begin{split}
\mu  \ &  \ \leq  \ \ts\frac{1}{\vunderline{\varsigma}}\,\E\Big(\ts\sup_{\theta\in\Theta}|\Gamma(\theta)|\Big)\\[0.2cm]
& \ \leq \  m^{4\beta-1}\E\|G\|_{F},
\end{split}
\end{equation*}
where the second step uses the norm inequalities $\|\cdot\|_{_{\scaleto{\!\triangle}{4pt}}}\leq \|\cdot\|\op\leq \|\cdot\|_F$. To conclude, will show $\E\|G\|_F\lesssim 1$ as follows. Since $G$ has the same covariance operator as $Y_1Y_1\ttop$, we have
\begin{equation*}
\begin{split}
(\E\|G\|_F)^2 & \ \leq \ \E \|G\|_F^2\\[0.2cm]
& \ = \ \E \Big\|Y_1Y_1\ttop \, - \, \E(Y_1Y_1\ttop)\Big\|_F^2\\[0.2cm]
& \ \lesssim \ \E \|Y_1\|_2^4 \\[0.2cm]
& \ = \ \bigg\| \sum_{j=1}^k \|\langle X_1,\phi_j\rangle\|_{L^2}^2 \,\zeta_{1j}^2\bigg\|_{L^2}^2\\[0.2cm]
&\ \leq \ \bigg(\sum_{j=1}^k \|\langle X_1,\phi_j\rangle\|_{L^2}^2\|\zeta_{1j}\|_{L^4}^2\bigg)^2 \\[0.2cm]
& \ \lesssim \  \ \Big(\ts\sum_{j=1}^k \lambda_j(\Sigma)\Big)^2 \\[0.2cm]
& \ \lesssim  \ 1.
\end{split}
\end{equation*}
Applying the last several bounds to~\eqref{eqn:anticintermed} yields the stated result. \qed

\begin{lemma}\label{lem:minvar}
Suppose the conditions of Theorem~\ref{thm:main} hold, and let $\underline\varsigma$ and $\bar\varsigma$ be as defined in~\eqref{eqn:varsigupperlower}. Then,
\begin{equation*}
\underline{\varsigma} \ \gtrsim \ m^{-4\beta+1} \ \ \ \ \ \ \text{and } \ \ \ \ \ \ \bar{\varsigma} \ \lesssim \ 1.
\end{equation*}
\end{lemma}
\proof 
Recall that the elements of the set $\Theta=\mathbb{B}_{\scaleto{\!\triangle}{4pt}}^k\times \{\pm 1\}$ are written as $\theta=(v,s)$, and that the vector $\zeta_1\in\R^k$ is defined in~\eqref{eqn:zetadef}. Then,
\begin{align*}
\varsigma^2(\theta)  & \ = \ \var(s  v\ttop Gv)\\[0.2cm]
& \ = \ \var\Big(\llangle Y_1Y_1\ttop, \ vv\ttop\rrangle\Big)\\[0.2cm]
& \ = \ \var\Big(\llangle \zeta_1\zeta_1\ttop,\, \Sigma_k^{1/2}vv\ttop  \Sigma_k^{1/2}\rrangle\Big)\\[0.2cm]
& \ = \  \Big\llangle \Sigma_k^{1/2}vv\ttop\Sigma_k^{1/2} \,\J_k(\Sigma_k^{1/2}vv\ttop\Sigma_k^{1/2})\Big\rrangle\\[0.2cm]
& \ \asymp \ \Big\|\Sigma_k^{1/2}vv\ttop\Sigma_k^{1/2}\Big\|_F^2 \text{  \ \ \ \ \ \ \ \ (Assumption~\ref{A:model}(c))}\\[0.2cm]
& \ = \ \langle v, \Sigma_k v\rangle^2.
\end{align*}

On one hand, the previous calculation of $\varsigma^2(\theta)$ implies $\bar\varsigma\lesssim 1$, since $\|v\|_2\leq 1$ and $\|\Sigma_k\|\op=\|\Sigma\|\op \lesssim 1$. On the other hand, recall that the condition $\theta=(v,s)\in\Theta$ implies $\sum_{j=1}^{m}\langle v,e_j\rangle^2\geq m^{-2\beta+1}$, and that $\Sigma_k=\textsf{Diag}(\lambda_1(\Sigma),\dots,\lambda_k(\Sigma))$. Therefore, we have the lower bound
\begin{align}
\langle v, \Sigma_k v\rangle^2 
&  \ \geq \ \bigg(\sum_{j=1}^{m} \lambda_{j}(\Sigma)\langle v,e_j\rangle^2\bigg)^2\\[0.2cm]
& \ \geq \ \Big(\lambda_{m}(\Sigma)m^{-2\beta+1}\Big)^2\\[0.2cm]
& \ \gtrsim \ m^{-8\beta+2}.
\end{align}
Hence, the lower bound $\varsigma(\theta)\gtrsim m^{-4\beta+1}$ holds uniformly over $\theta\in\Theta$, which completes the proof.\qed
~\\

\subsection{Approximating $\|\cdot\|\op$ with $\|\cdot\|_{_{\scaleto{\!\triangle}{4pt}}}$}\label{sec:II}
\begin{proposition}\label{prop:optri}
Suppose the conditions of Theorem~\ref{thm:main} hold. Then,
\begin{equation}\label{eqn:optri}
d_{\textup{K}}\Big(\mathcal{L}\big(\|G\|\op) \ , \ \mathcal{L}\big(\|G\|_{_{\scaleto{\!\triangle}{4pt}}}\big) \Big)
  \ \lesssim \ \ts\frac{1}{n}.
\end{equation}
\end{proposition}

\noindent\textbf{Remark.}\label{specialremark}
For the arguments used to prove Proposition~\ref{prop:optri} in this section, we may assume without loss of generality that $m<\dim(\mathbb{H})$, due to the following reasoning: Based on the definitions of $m$ and $k$ in~\eqref{eqn:mdef} and~\eqref{eqn:kdef}, if the condition $m<\dim(\mathbb{H})$ does not hold, then $m=k=\dim(\mathbb{H})$, and in this case the set $\mathbb{B}_{\scaleto{\!\triangle}{4pt}}^k$ defined in~\eqref{eqn:Bktriangledef} contains the unit sphere of $\R^k$. This implies that $\|\cdot\|\op$ and $\|\cdot\|_{\scaleto{\!\triangle}{4pt}}$ are equal on $\mathbb{S}^{k\times k}$, which means that the left side of~\eqref{eqn:optri} is 0.

\proof  Here, we only explain the proof at a high level, deferring the substantial steps to lemmas that are given later in this section. 
For any $t\in\R$, define the events
\begin{equation}\label{eqn:ABdefs}
\mathcal{A}(t)=\big\{\|G\|_{_{\scaleto{\!\triangle}{4pt}}} \leq t\big\} \text{ \ \ \ and \ \ \ } 
\mathcal{B}(t)=\big\{\|G\|_{\scaleto{\!\triangledown}{4pt}}> t\big\}.
\end{equation}
In terms of this notation, it is straightforward to check that the following relation holds for any $t\in\R$,
\begin{equation*}
\bigg |\P\big( \|G\|\op \leq t\big) \ - \  \P\big(\|G\|_{_{\scaleto{\!\triangle}{4pt}}}\leq t\big)\bigg |  \ = \ \P\Big(\mathcal{A}(t)\cap \mathcal{B}(t)\Big).
\end{equation*}
For any given pair of real numbers $t_1$ and $t_2$ satisfying $t_1\leq t_2$, it is possible to show that the following inclusion holds simultaneously for all $t\in\R$,
\begin{equation}\label{eqn:keyinclusion}
(\mathcal{A}(t)\cap \mathcal{B}(t)) \ \subset \ (\mathcal{A}(t_2)\cup \mathcal{B}(t_1)).
\end{equation}
This can be checked by separately considering the cases $t>t_1$ and $t\leq t_1$. In the case when $t>t_1$, we have $\mathcal{B}(t)\subset\mathcal{B}(t_1)$, and hence
$$
\big(\mathcal{A}(t ) \cap \mathcal{B}(t)\big) \subset \mathcal{B}(t)\subset \mathcal{B}(t_1)\subset \big(\mathcal{A}(t_2)\cup\mathcal{B}(t_1)\big).$$
Secondly, when $t\leq t_1$, we must also have $t\leq t_2$, which implies $\mathcal{A}(t)\subset\mathcal{A}(t_2)$, and then analogous reasoning can be used.
After applying a union bound to~\eqref{eqn:keyinclusion} and taking a supremum over $t\in\R$, we have
$$d_{\textup{K}}\big(\mathcal{L}\big(\| G\|\op\big) \ , \ \mathcal{L}\big(\|G\|_{_{\scaleto{\!\triangle}{4pt}}}\big) \big)
 \ \leq \ \P(\mathcal{A}(t_2))\,+\, \P(\mathcal{B}(t_1)).
 $$
The lower-tail probability $\P(\mathcal{A}(t_2))$ and the upper-tail probability $\P(\mathcal{B}(t_1))$ will be analyzed respectively in Lemmas~\ref{lem:lowertail} and~\ref{lem:uppertail} below, and suitable values of $t_1$ and $t_2$ will be specified in those results so that~\eqref{eqn:optri} holds.  Specifically, the values will be chosen so that $t_1\asymp m^{-2\beta+1}\sqrt{\log(n)}$ and $t_2\asymp \ell^{-2\beta}\sqrt{\log(\lceil \ell^{1/2} \rceil)}$. To check that these values satisfy the constraint $t_1\leq t_2$ for all large $n$,
 first note that the remark given before the current proof allows us assume that $m<\dim(\mathbb{H})$, which implies that $m$ is at least a fractional power of $n$. Consequently the definitions of $m$ and $\ell$ in~\eqref{eqn:mdef} and~\eqref{eqn:elldef} imply $t_1/t_2\to 0$ as $n\to\infty$.\qed

~\\

\subsubsection{Upper bound for $\P(\mathcal{A}(t_2))$}
To upper bound the probability $\P(\mathcal{A}(t_2))$ in the proof of Proposition~\ref{prop:optri}, we will need a specialized lower-tail bound for the maxima of correlated Gaussian random variables, which was previously developed in Theorem 2.3 of~\citep{Lopes:EJS}. For this purpose, we define the \emph{stable rank} of any non-zero positive semidefinite matrix $A$ as 
\begin{equation*}
{\tt{r}}(A) \ = \ \ts\frac{\tr(A)^2}{\|A\|_F^2}.
\end{equation*}
To allow the lower-tail bound to be understood more clearly, we have stated it below in terms of an index $N$ that is distinct from the sample size $n$ in the current paper.

\begin{lemma}[\cite{Lopes:EJS} Theorem 2.3]\label{lem:key}
Let $(\gamma_1,\dots,\gamma_{N})$ be a Gaussian random vector with a correlation matrix $R$, and entries satisfying $\E(\gamma_j)=0$ and $\var(\gamma_j)=1$ for all $j=1,\dots,N$. 
Also let $a,b \in (0,1)$ be any constants that are fixed with respect to $N$. Then, there is a constant $C>0$ depending only on $a$ and $b$ such that the following inequality holds for any integer $l_N$ satisfying $2\leq l_N\leq \ts\frac{a^2}{4} {{\tt{r}}}(R)$,
\begin{equation}\label{eqn:mainrefined}
\P\bigg(\max_{1\leq j\leq N} \gamma_j \, \leq \,  b \sqrt{2(1-a)\log(l_N)}\bigg) \ \leq \ C\, l_N^{\frac{-(1-a)(1-b)^2}{a}}\log(l_N)^{\frac{1-a(2-b)-b}{2a}}.
\end{equation}
\end{lemma}
~\\
\noindent We now state the main result of this subsection, which yields the desired upper bound on $\P(\mathcal{A}(t_2))$.

\begin{lemma}\label{lem:lowertail} 
Suppose the conditions of Theorem~\ref{thm:main} hold, and that $m<\dim(\mathbb{H})$. Then, there is a constant $c_2>0$ not depending on $n$ such that the choice
\begin{equation}\label{eqn:t2def}
t_2=c_2\ell^{-2\beta}\ts\sqrt{\log(\lceil \ell^{1/2}\rceil)}
\end{equation}
implies
\begin{equation*}
\P(\mathcal{A}(t_2)) \ \lesssim \ \ts\frac{1}{n}.
\end{equation*}
\end{lemma}
\proof  Based on the definition of $\|\cdot\|_{_{\scaleto{\!\triangle}{4pt}}}$ in~\eqref{eqn:triupdef}, the inequality $ \|A\|_{_{\scaleto{\!\triangle}{4pt}}}\geq \max_{1\leq j\leq \ell}\llangle A, e_je_j\ttop \rrangle $ holds for any $A\in\mathbb{S}^{k\times k}$, and so 
\begin{align*}
\small
\P(\mathcal{A}(t_2))  
& \ \leq \  \P\Big(\max_{1\leq j\leq \ell}\llangle G, e_je_j\ttop \rrangle \leq t_2\Big).
\end{align*}
(The choice to take the maximum over $j=1,\dots,\ell$ rather than $j=1,\dots,m$ is made so that the magnitudes of $t_2$ and $t_1$ in~\eqref{eqn:t2def} and~\eqref{eqn:t1def} will satisfy $t_1\leq t_2$ when $n$ is large, as needed in the proof of Proposition~\ref{prop:optri}.) 
Let $g_j=\llangle G,e_je_j\ttop\rrangle$ for $j=1,\dots,\ell$, and $\tau = \min_{1\leq j\leq \ell}\|g_j\|_{L^2}$. Also define the standardized variables $\bar g_j=g_j/\|g_j\|_{L^2}$. Then, it is straightforward to check that
\begin{equation}\label{eqn:barg}
\P\Big(\max_{1\leq j\leq \ell}\llangle G, e_je_j\ttop \rrangle \leq t_2\Big) \ \leq \ \P\Big(\max_{1\leq j\leq \ell} \bar g_j  \leq \ts\frac{t_2}{\tau}\Big).
\end{equation}
To control $t_2/\tau$, Lemma~\ref{lem:Rfrob} will show that $\tau \ \gtrsim \ \ell^{-2\beta}$. Hence, for any numbers $a,b\in (0,1)$ that are fixed with respect to $n$, there is a constant $c_2>0$ that may be used in the definition  of $t_2$ in \eqref{eqn:t2def} such that 
\begin{equation}\label{eqn:abbound}
\ts\frac{t_2}{\tau} \ \leq \ b\sqrt{2(1-a)\log(\lceil \ell^{1/2}\rceil)}.
\end{equation}
We are almost in position to apply the lower-tail bound from Lemma~\ref{lem:lowertail} to the right side of~\eqref{eqn:barg}, with $\bar g_j$ playing the role of $\gamma_j$,  $\ell$ playing the role of $N$, and $\lceil \ell^{1/2}\rceil $ playing the role of $l_N$. One more ingredient we need is a lower bound on the stable rank of the correlation matrix $R$ of the Gaussian vector $(\bar g_1,\dots,\bar g_{\ell})$. Specifically, Lemma~\ref{lem:Rfrob} below shows that
\begin{equation}\label{eqn:ranklowertemp}
{\tt{r}}(R) \ \gtrsim \ \ell.
\end{equation}
By the assumption that $m<\dim(\mathbb{H})$, it follows that $\ell$ is lower bounded by a fractional power of $n$. Consequently, we have $\lceil \ell^{1/2}\rceil/\ell\to 0$ as $n\to\infty$. Moreover,~\eqref{eqn:ranklowertemp} implies that for any fixed $a\in(0,1)$, we have $2\leq \lceil \ell^{1/2}\rceil \leq \frac{a^2}{4}{\tt{r}}(R)$ for all large $n$. So, combining Lemma~\ref{lem:key} with the bound~\eqref{eqn:abbound} gives
\begin{equation}\label{eqn:almost}
\P\Big(\max_{1\leq j\leq \ell} \bar g_j  \leq \ts\frac{t_2}{\tau}\Big)
  \ \lesssim \ \lceil \ell^{1/2}\rceil^{-\frac{(1-a)(1-b)^2}{a}}\log(\lceil \ell^{1/2}\rceil)^{\frac{1-a(2-b)-b}{2a}}.
\end{equation}
Finally, since the number $\ell$ is at least a fractional power of $n$, and since $a$ can be chosen arbitrarily small, the right side of~\eqref{eqn:almost} can be made to be at most of order $1/n$.\qed

\begin{lemma}\label{lem:Rfrob}
Suppose the conditions of Theorem~\ref{thm:main} hold. Then, 
\begin{equation}\label{eqn:varlower}
\min_{1\leq j\leq \ell} \var\Big(\llangle G, e_je_j\ttop \rrangle\Big) \ \gtrsim \ \lambda_{\ell}(\Sigma)^2.
\end{equation}
In addition, if $R$ denotes the correlation matrix of the random vector $(\llangle G, e_1e_1\ttop\rrangle,\dots,\llangle G, e_{\ell}e_{\ell}\ttop\rrangle)$, then
\begin{equation}\label{eqn:corlower}
{\tt{r}}(R) \ \gtrsim \  \ell.
\end{equation}

\end{lemma}
\proof 
First we establish~\eqref{eqn:varlower}. Recall the notation $\zeta_1=(\zeta_{11},\dots,\zeta_{1k})$, as well as the facts that $e_j\ttop Y_1=\sqrt{\lambda_j(\Sigma)}e_j\ttop\zeta_1$, and that $G$ has the same covariance operator as $Y_1Y_1\ttop$. So, for any $j=1,\dots,\ell$, we have
\begin{equation}\label{eqn:varlowerderive}
\begin{split}
\var(\llangle G,e_je_j\ttop\rrangle) 
& \ = \ \var(\llangle Y_1Y_1\ttop, e_je_j\ttop\rrangle)\\
& \ = \ \lambda_j(\Sigma)^2\var\big(\llangle \zeta_1\zeta_1\ttop\,,\, e_je_j\ttop\rrangle\big)\\
& \ = \ \lambda_j(\Sigma)^2 \big\llangle  e_je_j\ttop,  \J_k (e_je_j\ttop)\big\rrangle\\
& \ \gtrsim \ \lambda_{j}(\Sigma)^2,
\end{split}
\end{equation}
where the last step follows from Assumption~\ref{A:model}(c). This establishes~\eqref{eqn:varlower}.

Regarding~\eqref{eqn:corlower}, it suffices to show that $\|R\|_F^2\lesssim \ell$, since $\tr(R)=\ell$.
To proceed, recall the notation $g_j=\llangle G,e_je_j\ttop\rrangle$, and note that the squared entries of $R$ are given by
\begin{equation}\label{eqn:corformula}
R_{jj'}^2 \ = \ \frac{\cov(g_j,g_{j'})^2}{\var(g_j)\var(g_{j'})}.
\end{equation}
By essentially repeating the first several steps of~\eqref{eqn:varlowerderive}, we have
\begin{equation*}
\begin{split}
\cov(g_j,g_{j'}) 
& \ = \ \lambda_j(\Sigma)\lambda_{j'}(\Sigma) \big\llangle  e_je_j\ttop,  \J_k ( e_{j'}e_{j'}\ttop)\big\rrangle.
\end{split}
\end{equation*}
If we let $K\in\mathbb{S}^{\ell \times \ell}$ have entries  $K_{jj'} = \big \llangle  e_je_j\ttop,  \J_k ( e_{j'}e_{j'}\ttop)\big\rrangle$, then~\eqref{eqn:varlowerderive} and~\eqref{eqn:corformula} imply that the squared entries of $R$ satisfy the bound
$R_{jj'}^2 \ \lesssim \ K_{jj'}^2$.
This yields
\begin{equation}\label{eqn:Kholder}
\begin{split}
\|R\|_F^2 & \ \lesssim \ \|K\|_F^2  \ \leq \ \ell\, \|K\|\op^2.
\end{split}
\end{equation}
Next, observe that if the operator $\J_k$ is represented as a symmetric matrix of size $\frac{1}{2}k(k+1)\times \frac{1}{2}k(k+1)$ with entries $\{\llangle  e_ie_{j}\ttop,  \J_k ( e_{i'}e_{j'}\ttop)\rrangle\}$ for $ i\leq j$ and $i'\leq j'$, then $K$ corresponds to a  principal submatrix of $\mathsf{J}_k$. Therefore, $\|K\|\op\leq \|\mathsf{J}_k\|\op\lesssim 1$ by Assumption~\ref{A:model}(c), and so $\|R\|_F^2\lesssim \ell$.\qed

\subsubsection{Upper bound for \normalfont{$\P(\mathcal{B}(t_1))$}}

\begin{lemma}\label{lem:uppertail}Suppose the conditions of Theorem~\ref{thm:main} hold. Then, there is a constant $c_1>0$ not depending on $n$, such that the choice 
\begin{equation}\label{eqn:t1def}
t_1=c_1m^{-2\beta+1}\sqrt{\log(n)}
\end{equation}
implies
\begin{equation}\label{eqn:PBt1}
\P(\mathcal{B}(t_1)) \ \lesssim \ \ts\frac{1}{n}.
\end{equation}
\end{lemma}
\proof The overall idea for proving \eqref{eqn:PBt1} is to develop an upper bound on the $L^q$ norm of the random variable \smash{$\|G\|_{\scaleto{\!\triangledown}{4pt}}$} and apply Chebyshev's inequality.
To begin with the details, let $\eta=m^{-\beta+1/2}$, and define the diagonal linear operator $D:\R^k\to\R^k$ according to 
\begin{equation}\label{eqn:Ldef}
D x \ = \ \sum_{j=1}^{m} 2\eta \langle x,e_j\rangle e_j \ + \  \sum_{l= m+1}^k 2\langle x,e_l\rangle e_l.
\end{equation}
In a moment, we will need to use the fact that the image of $\mathbb{B}_{\scaleto{\!\triangledown}{4pt}}^k$ under the inverse $D^{-1}$ is contained in the unit ball $\mathbb{B}^k$. This can be seen by noting that for any $u\in \mathbb{B}_{\scaleto{\!\triangledown}{4pt}}^k$ we have
\begin{equation*}
\begin{split}
\|D^{-1}u\|_2 & \ = \  \bigg\|\ts\frac{1}{2\eta}\sum_{j=1}^{m} \langle u,e_j\rangle e_j +\ts\frac{1}{2}\sum_{l= m+1}^k\langle u,e_l\rangle e_l\bigg\|_2\\[0.2cm]
& \ \leq \ \ts\frac{1}{2\eta} \Big(\sum_{j=1}^{m} \langle u,e_j\rangle^2\Big)^{1/2} + \ts\frac{1}{2}\|u\|_2\\[0.2cm]
& \ \leq \ 1,
\end{split}
\end{equation*}
where the last step uses the definitions of $\mathbb{B}_{\scaleto{\!\triangledown}{4pt}}^k$ and $\eta$.
Hence, for any $u\in\mathbb{B}_{\scaleto{\!\triangledown}{4pt}}^k$, there is a vector $w\in \mathbb{B}^k$ such that $Dw=u$, which yields
\begin{equation*}
\begin{split}
\|G\|_{\scaleto{\!\triangledown}{4pt}}
& \ \leq  \ \sup_{w\in \mathbb{B}^k} \Big| \Big\langle Dw,  G Dw\Big\rangle\Big| \\[0.2cm]
& \ \leq \ \| DGD\|_F,
\end{split}
\end{equation*}
where for ease of notation we have identified $D$ with its matrix representation in the standard basis.

We are now in position to bound the $L^q$ norm of $\|G\|_{\scaleto{\!\triangledown}{4pt}}$, taking $q=\log(n)\vee 3$. Since the $(i,j)$ entry of $DGD$ is equal to $G_{ij}D_{ii}D_{jj}$, we have
\begin{equation}\label{eqn:firstlqgtri}
\begin{split}
\big\|\|G\|_{\scaleto{\!\triangledown}{4pt}}\big\|_{L^q} 
& \ \leq \ \bigg(\sum_{i=1}^k\sum_{j=1}^k \|G_{ij}\|_{L^q}^{2} D_{ii}^2 D_{jj}^2\bigg)^{1/2}.
\end{split}
\end{equation}
Using a standard moment bound for centered Gaussian random variables~\cite[p.25]{Vershynin:2018}, as well as the fact that $G$ has the same covariance operator as $Y_1Y_1\ttop$,  it follows that 
\begin{equation*}
\begin{split}
 \|G_{ij}\|_{L^q}^2 & \ \lesssim \ q \var(G_{ij})\\
 & \ = \ q \var(Y_{1i}Y_{1j})\\
 & \ \leq \ q \,\lambda_i(\Sigma)\lambda_j(\Sigma)\E(\zeta_{1i}^2\zeta_{1j}^2)\\
 &\ \lesssim  \ q \,\lambda_i(\Sigma)\lambda_j(\Sigma),
\end{split}
\end{equation*}
where we recall that $Y_1=(\sqrt{\lambda_1(\Sigma)}\zeta_{11},\dots,\sqrt{\lambda_k(\Sigma)}\zeta_{1k})$. Combining with~\eqref{eqn:firstlqgtri}, we conclude that
\begin{equation*}
\begin{split}
\big\|\|G\|_{\scaleto{\!\triangledown}{4pt}}\big\|_{L^q}
& \ \lesssim \ \sqrt{q}\sum_{j=1}^k \lambda_j(\Sigma) D_{jj}^2\\
& \ \lesssim \  \sqrt q\bigg(\displaystyle \sum_{j=1}^m \eta^2 j^{-2\beta}  + \sum_{j=m+1}^k  j^{-2\beta}\bigg)\\
& \ \lesssim \ \sqrt q\big(\eta^2  \ + \ m^{-2\beta+1}\big)\\
& \ \lesssim \sqrt{\log(n)} m^{-2\beta+1}.
\end{split}
\end{equation*}
Under the choice of $q=\log(n)\vee 3$, this bound on $\big\|\|G\|_{\scaleto{\!\triangledown}{4pt}}\big\|_{L^q}$ leads to the stated result via Chebyshev's inequality for a suitable choice of $c_1$.\qed

\section{The terms \protect{$\II_n$ and $\hat{\II}_n$}: Gaussian and bootstrap approximation}\label{sec:III}
\begin{proposition}\label{prop:optriI}
Suppose the conditions of Theorem~\ref{thm:main} hold. Then, 
\begin{equation}\label{eqn:propIIIbound}
\II_n \ \lesssim  \ n^{-\frac{\beta-1/2}{2\beta+4+\e}}.\tag{i}
\end{equation}
In addition, there is a constant $c>0$ not depending on $n$ such that the event
\begin{equation}\label{eqn:propIIIhatbound}
\hat{\II}_n \ \leq  \ c\,n^{-\frac{\beta-1/2}{2\beta+4+\e}}\tag{ii}
\end{equation}
holds with probability at least $1-c/n$.
\end{proposition}
\proof We begin with part (i).  Let $\mathscr{C}$ denote the class of all convex Borel subsets of $\mathbb{S}^{k\times k}$. Since the function $\|\cdot\|\op$ is convex, we have
\begin{equation}\label{eqn:boundIIIfirst}
\begin{split}
\II_n  & \ = \ \sup_{t\in\R}\bigg|\P\Big(\sqrt n\|\hat\Sigma_k-\Sigma_k\|\op \leq t\Big) \ - \ \P\big(\|G\|\op \leq t\big)\bigg|\\[0.3cm]
& \ \leq \ \sup_{\mathcal{C}\in\mathscr{C}} \bigg|\P\Big(\sqrt n(\hat\Sigma_k-\Sigma_k) \in\mathcal{C} \Big) \ - \ \P\big(G \in\mathcal{C}\big)\bigg|.
\end{split}
\end{equation}
Next, recall that $\hat\Sigma_k-\Sigma_k = \ts\frac{1}{n}\sum_{i=1}^n Y_iY_i\ttop-\E(Y_iY_i\ttop)$, and let ${\mathsf{C}}:\mathbb{S}^{k\times k}\to\mathbb{S}^{k\times k}$ denote the covariance operator of $Y_1Y_1\ttop$, so that for any $A,B\in\mathbb{S}^{k\times k}$ we have
$$\llangle A, {\mathsf{C}}(B)\rrangle  \ = \ \cov\Big(\big\llangle A,Y_1Y_1\ttop \big\rrangle, \big\llangle B, Y_1Y_1\ttop\big\rrangle\Big).$$
It is possible to write the operator $\mathsf{C}$ more explicitly as
\begin{equation}\label{eqn:Cinvertibility}
\mathsf C = \big(\Sigma_k^{1/2}\ok \Sigma_k^{1/2}\big)\,\J_k\, \big(\Sigma_k^{1/2}\ok \Sigma_k^{1/2}\big),
\end{equation}
which can be verified by observing that the identity $Y_k=\Sigma_k^{1/2}\zeta_1$ and the definition of $\mathsf{J}_k$ yield
\begin{equation*}
\begin{split}
\Big\llangle A\, ,\, \big(\Sigma_k^{1/2}\ok \Sigma_k^{1/2}\big)\,\J_k\, \big(\Sigma_k^{1/2}\ok \Sigma_k^{1/2}\big)(B)\Big\rrangle \ & = \Big\llangle \Sigma_k^{1/2} A\Sigma_k^{1/2}, \mathsf{J}_k\big(\Sigma_k^{1/2} B\Sigma_k^{1/2}\big)\Big\rrangle\\[0.2cm]
& \ = \ \cov\Big(\big\llangle \Sigma_k^{1/2} A\Sigma_k^{1/2}, \zeta_1\zeta_1\ttop\big\rrangle \ , \ \big\llangle \Sigma_k^{1/2} B\Sigma_k^{1/2}, \zeta_1\zeta_1\ttop\big\rrangle\Big)\\[0.2cm]
& \ = \ \cov\Big(\big\llangle A,Y_1Y_1\ttop \big\rrangle, \big\llangle B, Y_1Y_1\ttop\big\rrangle\Big).
\end{split}
\end{equation*}

The representation of $\mathsf{C}$ in~\eqref{eqn:Cinvertibility} shows that it is a composition of invertible operators, and hence invertible.
So, for each $i=1,\dots,n$, we may define 
\begin{equation}\label{eqn:Midef}
M_i={\mathsf{C}}^{-1/2}\Big(Y_iY_i\ttop-\E(Y_iY_i\ttop)\Big),
\end{equation}
which we will regard as a centered and isotropic random vector in $\mathbb{S}^{k\times k}$. (Note that $M_i$ takes values in~$\mathbb{S}^{k\times k}$ since $\mathsf{C}^{-1/2}$ is an operator from $\mathbb{S}^{k\times k}$ to itself.)  Using this definition, the bound~\eqref{eqn:boundIIIfirst} leads to
\begin{equation}\label{eqn:boundIII2nd}
\begin{split}
\II_n   \ \leq \  \sup_{\mathcal{C}\in\mathscr{C}} \bigg| \P\Big(\ts\frac{1}{\sqrt{n}}(M_1+\dots+M_n)\in \mathcal{C}\Big) \ - \ \gamma(\mathcal{C})\bigg|,
\end{split}
\end{equation}
where $\gamma$ denotes the standard Gaussian distribution on $\mathbb{S}^{k\times k}$.  To bound the last expression, we will apply a version of the multivariate Berry-Esseen Theorem due to~\citep{Fang:Koike:Balls}, recorded later in Lemma~\ref{lem:fangkoike}, which yields
\begin{equation}\label{eqn:FK}
\begin{split}
\II_n  & \ \lesssim \ \displaystyle \ts\frac{1}{n^{1/2}} \cdot \Big(\ts\frac{k(k+1)}{2}\Big)^{1/4}\cdot (\E\|M_1\|_F^4)^{1/2} \cdot \big(|\log(\E\|M_1\|_F^4/n)|\vee 1\big).
\end{split}
\end{equation}
Finally, in Lemma~\ref{lem:F4bound},  it is shown that
$\E\|M_1\|_F^4 \ \lesssim \ k^4$, and combining this with $k\lesssim n^{\frac{1}{2\beta+4}}$ leads to the bound (i).\\

Turning to part (ii), recall that $Y_1^*,\dots,Y_n^*\in\R^k$ are sampled with replacement from $Y_1,\dots,Y_n$. In addition, let $\hat{{\mathsf{C}}}:\mathbb{S}^{k\times k}\to\mathbb{S}^{k\times k}$ denote the covariance operator of the conditional distribution of $Y_1^*(Y_1^*)\ttop$ given $X_1,\dots,X_n$.
 It follows from Lemma~\ref{lem:Chatinverse} that $\hat{\mathsf{C}}$ is invertible with probability at least $1-c/n$, and on this high probability event we may define the bootstrap counterpart of $M_i$ as
\begin{equation}\label{eqn:Mistardef}
M_i^*=\hat{{\mathsf{C}}}^{-1/2}\Big(Y_i^*(Y_i^*)\ttop \, - \, \E\big(Y_i^*(Y_i^*)\ttop \big|X\big)\Big).
\end{equation}
Next, by repeating the steps in the proof of part (i), there is a constant $c>0$ not depending on $n$ such that the following bound holds with probability at least $1-c/n$,
  \begin{equation*}
\begin{split}
\hat{\II}_n  
& \ \ \leq \ \ \displaystyle \ts\frac{c}{n^{1/2}}\cdot  \Big(\frac{k(k+1)}{2}\Big)^{1/4}\cdot (\E(\|M_1^*\|_F^4|X))^{1/2}\cdot\big(|\log(\E(\|M_1^*\|_F^4|X)/n)|\vee 1\big).
\end{split}
\end{equation*}
\normalsize
Lastly, it is shown in Lemma~\ref{lem:conditionalFnorm} that the bound
 $$\E(\|M_1^*\|_F^4|X) \ \leq \ c\,k^4 \log(n)^{8\nu}$$
 holds with probability at least $1-c/n$, which leads to the statement (ii).
\qed

~\\

\begin{lemma}\label{lem:F4bound} If the conditions of Theorem~\ref{thm:main} hold and $M_1$ is defined as in~\eqref{eqn:Midef}, then
$$\E\big\|M_1\|_F^4 \ \lesssim \ k^4.$$
\end{lemma}

\proof  Due to~\eqref{eqn:Cinvertibility}, the inverse of $\mathsf{C}$ is given by
\begin{equation*}
\mathsf C^{-1} = \big(\Sigma_k^{-1/2}\ok \Sigma_k^{-1/2}\big)\,\J_k^{-1}\, \big(\Sigma_k^{-1/2}\ok \Sigma_k^{-1/2}\big).
\end{equation*}
Recall the notation $\zeta_i=(\zeta_{i1},\dots,\zeta_{ik})$ for each $i=1,\dots,n$, and that this random vector satisfies $\zeta_i=\Sigma_k^{-1/2}Y_i$.
Using Assumption~\ref{A:model}(c), the following bounds hold with probability 1,
\begin{equation}\label{eqn:Frobinter}
\begin{split}
\|M_1\|_F^4 
& \ = \ \Big\llangle Y_1Y_1\ttop -\E(Y_1Y_1\ttop) , {\mathsf{C}}^{-1}\Big(Y_1Y_1\ttop-\E(Y_1Y_1\ttop)\Big)\Big\rrangle^2\\[0.2cm]
& \ = \ \Big\llangle \zeta_1\zeta_1\ttop - I , {\J}_k^{-1}\big(\zeta_1\zeta_1\ttop - I\big)\Big\rrangle^2 \\[0.2cm]
& \ \leq \ \ts\frac{1}{\lambda_{\min}(\J_k)^2} \big\|\zeta_1\zeta_1\ttop - I\big\|_F^4\\[0.2cm]
& \ \leq \  c\Big(\|\zeta_1\|_2^8  \, + \, k^2\Big),
\end{split}
\end{equation}
for some constant $c>0$ not depending on $n$. To complete the proof, it remains to bound $\E\|\zeta_1\|_2^8$, for which we have
\begin{equation*}
\begin{split}
\Big(\E\|\zeta_1\|_2^8\Big)^{1/4} 
& \ \leq \ \sum_{j=1}^k\|\zeta_{1j}^2\|_{L^4} \ \lesssim \ k,
\end{split}
\end{equation*}
as needed.   \qed

\begin{lemma}\label{lem:conditionalFnorm} Suppose that the conditions of Theorem~\ref{thm:main} hold, and let $M_1^*$ be as defined in~\eqref{eqn:Mistardef}. Then, there is a constant $c>0$ not depending on $n$ such that the following event holds with probability at least $1-c/n$,
\begin{equation*}
\E\big(\|M_1^*\|_F^4\,\big|\,X\big)  \ \leq \ c\, k^4\log(n)^{8\nu}.
\end{equation*}
\end{lemma}

\proof The proof proceeds by analogy with that of Lemma~\ref{lem:F4bound}, but involves some extra technical hurdles. First, let $\zeta_1^*=\Sigma_k^{-1/2}Y_1^*$, and let $\hat{\!\J}_k$ denote the conditional covariance operator of $\zeta_1^*(\zeta_1^*)\ttop$ given $X_1,\dots,X_n$, which satisfies 
\begin{equation}\label{eqn:hatJdef}
\llangle A\, , \, \hat{\!\J}_k(B)\rrangle \ = \  \ \cov\Big(\big\llangle A,\zeta_1^*(\zeta_1^*)\ttop  \big\rrangle\, , \,  \big\llangle B, \zeta_1^*(\zeta_1^*)\ttop \big\rrangle\, \Big|\,X\Big)
\end{equation}
for any $A,B\in\mathbb{S}^{k\times k}$. Also, recall that $\hat{\mathsf{C}}$ denotes the covariance operator of $Y_1^*(Y_1^*)\ttop$ conditional on $X_1,\dots,X_n$,  which is related to $\hat{\!\J}_k$ through the formula
\begin{equation}\label{eqn:hatCformula}
\hat{\mathsf{C}}=(\Sigma_k^{1/2}\ok \Sigma_k^{1/2})\ \hat{\!\J}_k^{}\, (\Sigma_k^{1/2}\ok \Sigma_k^{1/2}).
\end{equation}
Hence, $\hat{\mathsf{C}}$ is invertible whenever $\hat{\!\J}_k$ is, and it will be shown in Lemma~\ref{lem:Chatinverse} that there is a constant $c>0$ not depending on $n$ such that the event $\{\lambda_{\min}(\,\hat{\!\J}_k\,)\geq 1/c\}$ holds with probability at least $1-c/n$. When this event occurs, we may bound $\|M_1^*\|_F^2$ according to 
\begin{equation}\label{eqn:initialMstar}
\begin{split}
\|M_1^*\|_F^4
& \ = \ \bigg\llangle Y_1^*(Y_1^*)\ttop-\E(Y_1^*(Y_1^*)\ttop|X)\, , \, \hat{\mathsf{C}}^{-1}\Big((Y_1^*(Y_1^*)\ttop-\E(Y_1^*(Y_1^*)\ttop|X)\Big)\bigg\rrangle^2\\[0.2cm]
& \ = \ \bigg\llangle \zeta_1^*(\zeta_1^*)\ttop-\E(\zeta_1^*(\zeta_1^*)\ttop|X) \, , \, \hat{\!\J}_k^{\ -1}\Big(\zeta_1^*(\zeta_1^*)\ttop-\E(\zeta_1^*(\zeta_1^*)\ttop|X)\Big)\bigg\rrangle^2\\[0.2cm]
& \ \leq \ c\Big\|\zeta_1^*(\zeta_1^*)\ttop -\E(\zeta_1^*(\zeta_1^*)\ttop|X)\Big\|_F^4\\[0.3cm]
& \ \leq  \ c\|\zeta_1^*\|_2^8 \ + \ c\E\big(\|\zeta_1^*\|_2^8\,\big|X\big).
\end{split}
\end{equation}
Therefore, the following event also occurs with probability at least $1-c/n$,
\begin{equation}\label{eqn:almost4th}
\begin{split}
\E\big(\|M_1^*\|_F^4\big|X\big) & \ \leq \ c\, \E\big(\|\zeta_1^*\|_2^8\,\big|X\big)\\[0.2cm]
& \ = \ \ts\frac{c}{n}\sum_{i=1}^n \|\zeta_i\|_2^8.
\end{split}
\end{equation}
So, to obtain a high probability bound on $\E\big(\|M_1^*\|_F^4\big|X\big)$ it suffices to obtain a high probability bound on $\ts\frac{1}{n}\sum_{i=1}^n \|\zeta_i\|_2^8$. For this purpose, let $q=\log(n)\vee 1$ and consider the moment bound
\begin{equation*}
\begin{split}
 \Big\|\ts\frac{1}{n}\sum_{i=1}^n \|\zeta_i\|_2^8\Big\|_{L^q} 
 & \ \leq \ \big\|\|\zeta_1\|_2^2\big\|_{L^{4q}}^4\\
   & \ \leq \ \bigg(\sum_{j=1}^k\|\zeta_{1j}\|_{L^{8q}}^2\bigg)^4\\
  & \ \lesssim \ q^{8\nu} k^4,
 \end{split}
\end{equation*}
where the last step uses Assumption~\ref{A:model}(a). Hence, by Chebyshev's inequality, the event $\frac{1}{n}\sum_{i=1}^n \|\zeta_i\|_2^8\geq cq^{8\nu}k^4$  occurs with probability at most $c/n$. Applying this to~\eqref{eqn:almost4th} completes the proof.\qed

\begin{lemma}\label{lem:Chatinverse}
Suppose the conditions of Theorem~\ref{thm:main} hold. Then, there is a constant $c>0$ not depending on $n$ such that the following events each hold with probability at least $1-c/n$,
\begin{equation}\label{eqn:hatJerror}
\|\,\hat{\!\J}_k-\J_k\|\op \ \leq \ \frac{  ck \log(n)^{2\nu+1/2}}{n^{1/2}},
\end{equation}
and
\begin{equation}\label{eqn:hatJlmin}
\lambda_{\min}(\,\hat{\!\J}_k) \ \geq \ 1/c.
\end{equation}
In addition, if the event~\eqref{eqn:hatJlmin} holds, then the operator $\hat{\mathsf{C}}$ given by the formula~\eqref{eqn:hatCformula} is invertible.
\end{lemma}
\proof Due to Weyl's inequality~\citep[Theorem 4.3.1]{Horn:Johnson}  and Assumption \ref{A:model}(c), the statement~\eqref{eqn:hatJlmin} follows from~\eqref{eqn:hatJerror}. Also, the assertion that~\eqref{eqn:hatJlmin} implies the invertibility of $\hat{\mathsf{C}}$ is immediate from the formula~\eqref{eqn:hatCformula}.

It remains to prove~\eqref{eqn:hatJerror}. Let $\overline{\zeta\zeta\ttop\!\!}\ =\frac{1}{n}\sum_{i=1}^n \zeta_i\zeta_i\ttop$, and note from~\eqref{eqn:hatJdef} that the operator $\hat{\!\J}_k$ can be expressed as
\begin{equation*}
\begin{split}
\hat{\!\J}_k & \ = \ \frac{1}{n}\sum_{i=1}^n (\zeta_i\zeta_i\ttop  -\overline{\zeta\zeta\ttop\!\!}\ )\otimes (\zeta_i\zeta_i\ttop  -\overline{\zeta\zeta\ttop\!\!} \ ).
\end{split}
\end{equation*}
To decompose $\hat{\!\J}_k$, define 
\begin{equation*}
\begin{split}
\hat{\!\J}_k\,\!\!'   \ = \   \ts\frac 1n\displaystyle \sum_{i=1}^n (\zeta_i\zeta_i\ttop - I)\otimes (\zeta_i\zeta_i\ttop - I) \ \ \ \ \ \text{ and } \ \ \ \ \ \ \ 
\hat{\!\J}_k\,\!\!'' & \ = \ (\overline{\zeta\zeta\ttop\!\!} - I)\otimes (\overline{\zeta\zeta\ttop\!\!}-I)
\end{split}
\end{equation*}
so that $\,\hat{\!\J}_k  \, =\, \hat{\!\J}_k\,\!\!' \, -\,  \hat{\!\J}_k\,\!\!'' $, which yields the bound
\begin{equation}\label{eqn:hatJsplitbound}
 \|\,\hat{\!\J}_k-\J_k\|\op \ \leq \ \|\,\hat{\!\J}_k\!\!'-\J_k\|\op \ + \ \|\,\hat{\!\J}_k\,\!\!''\|\op.
\end{equation}
For the first term on the right side,
note that if we put $\xi_i=\zeta_i\zeta_i\ttop-I$, then $\E(\xi_i\otimes \xi_i)=\J_k$, and so we may apply Lemma~\ref{lem:op} to obtain the following bound with $q=\log(n)\vee 3$,
\begin{equation*}
\begin{split}
\Big(\E\|\,\hat{\!\J}_k\,\!\!'-\J_k\|\op^q\Big)^{1/q}  & \ = \ \bigg(\E\Big\|\ts\frac{1}{n}\displaystyle \sum_{i=1}^n \xi_i\otimes \xi_i - \E(\xi_i\otimes \xi_i)\Big\|\op^q\bigg)^{1/q}\\[0.2cm]
& \ \lesssim \ \bigg( \ts\frac{\sqrt q}{\sqrt n}\, \|\J_k\|\op^{1/2}\, (\E\|\xi_1\|_F^{2q})^{\frac{1}{2q}}\bigg) \bigvee  \bigg(\ts\frac{q}{n}\, (\E \|\xi_1\|_F^{2q})^{\frac{1}{q}}\bigg).
\end{split}
\end{equation*}
Furthermore, we have
\begin{equation}\label{eqn:hatJinter}
\begin{split}
 (\E \|\xi_1\|_F^{2q})^{\frac{1}{2q}} 
 & \ \leq \  (\E \|\zeta_1\|_2^{4q})^{\frac{1}{2q}} \ + \ \sqrt k\\[0.2cm]
 & \ \leq \ \ts\sum_{j=1}^k \| \zeta_{1j}\|_{L^{4q}}^2\ + \ \sqrt k\\[0.2cm]
  & \ \lesssim \ q^{2\nu}k+ \ \sqrt k,
\end{split}
\end{equation}
where we have used Assumption~\ref{A:model}(a) in the last step. Combining this with $\|\J_k\|\op\lesssim1$ from Assumption~\ref{A:model}(c) yields
\begin{equation*}
\begin{split}
 \Big(\E\|\,\hat{\!\J}_k\,\!\!'-\J_k\|\op^q\Big)^{1/q} &  \ \lesssim \ \Big(\ts\frac{\sqrt q}{\sqrt n}(q^{2\nu}k+\sqrt k)\Big)\bigvee\Big(\ts\frac qn (q^{2\nu}k+\sqrt k)^2\Big)\\[0.2cm]
 & \ \lesssim  \ts\frac{q^{2\nu+1/2} k}{n^{1/2}}.
 \end{split}
\end{equation*}
To bound the second term on the right side of~\eqref{eqn:hatJsplitbound}, observe that
\begin{equation*}
\begin{split}
 \|\,\hat{\!\J}_k\,\!\!''\|\op & \ = \ \Big\|(\overline{\zeta\zeta\ttop\!\!} - I)\otimes (\overline{\zeta\zeta\ttop\!\!}-I)\Big\|\op\\[0.2cm]
 & \ = \ \|\overline{\zeta\zeta\ttop\!\!}-I\|_F^2\\[0.2cm]
 & \ \leq \ k \|\ts\frac{1}{n}\sum_{i=1}^n \zeta_i\zeta_i\ttop -I\|\op^2,
\end{split}
\end{equation*}
where we have used the fact that  $\|\cdot\|_F^2\leq k\|\cdot\|\op^2$ on $\mathbb{S}^{k\times k}$.
Since $\E(\zeta_1\zeta_1\ttop)=I$, we may apply Lemma~\ref{lem:op} to the last line, which yields
\begin{equation*}
\begin{split}
 \big(\E\|\,\hat{\!\J}_k\,\!\!''\|\op^q\big)^{1/q}
& \ \lesssim \ k\,\bigg[\Big( \ts\frac{\sqrt{2q}}{\sqrt n}\, (\E\|\zeta_1\|_2^{4q})^{\frac{1}{4q}}\Big) \bigvee  \Big(\ts\frac{2q}{n}\, (\E \|\zeta_1\|_2^{4q})^{\frac{1}{2q}}\Big)\bigg]^2\\[0.2cm]
& \ \lesssim  \ k \Big[\Big(\ts\frac{\sqrt{q}}{\sqrt n} \sqrt{q^{2\nu}k} \Big)\vee \Big(\ts\frac{q}{n}q^{2\nu}k\Big)\Big]^2 \\[0.2cm]
& \ \lesssim \ \ts\frac{q^{2\nu+1}k^2}{n},
\end{split}
\end{equation*}
where we have re-used part of~\eqref{eqn:hatJinter} in the second step.
Combining the last several steps and noting that $q^{2\nu+1}k^2/n\lesssim q^{2\nu+1/2}k/n^{1/2}$ gives 
\begin{equation*}
 \Big(\E \|\,\hat{\!\J}_k-\J_k\|\op^q\Big)^{1/q} \ \lesssim \ \ts\frac{q^{2\nu+1/2} k}{n^{1/2}}.
\end{equation*}
Therefore, the stated bound~\eqref{eqn:hatJerror} follows from Chebyshev's inequality. \qed

\section{The term $\hat{\III}_n$: Gaussian comparison}\label{sec:hatIV}

\begin{proposition}\label{prop:IVhat}
Suppose the conditions of Theorem~\ref{thm:main} hold. Then, there is a constant $c>0$ not depending on $n$ such that the event 
\begin{equation}\label{eqn:propIVstatement}
\hat{\III}_n \ \leq \  c\,n^{-\frac{\beta-1/2}{2\beta+4+\e}}
\end{equation}
holds with probability at least $1-c/n$.
\end{proposition}
\proof Define the operator $\hat{{\mathsf{D}}}:\mathbb{S}^{k\times k}\to\mathbb{S}^{k\times k}$, according to $\hat{\mathsf{D}}=\J_k^{-1/2}\,\,\hat{\!\J}_k\,\J_k^{-1/2}- \mathsf{I}$, where ${\mathsf{I}}$ is the identity operator.
As an initial step, we claim that the following bound holds with probability 1,
\begin{equation}\label{eqn:ivclaim}
\hat{\III}_n \ \leq \  2k\,\|\hat{{\mathsf{D}}}\|\op.
\end{equation}
In verifying the claim, we may restrict our attention to the case when the event $\{\|\hat{{\mathsf{D}}}\|\op\leq 1/2\}$ holds, since $\hat{\III}_n$ cannot be greater than 1. 
To proceed, we will adapt an argument from the proof of Lemma A.7 in~\citep{Spokoiny:2015}. Let $\mathscr{B}$ denote the collection of all Borel subsets of $\mathbb{S}^{k\times k}$, and let $d_{\text{KL}}$ denote the KL-divergence. Pinsker's inequality implies
\begin{equation*}
\begin{split}
\hat{\III}_n & \ = \ \sup_{t\in\R}\bigg|\P\big(\|G\|\op \leq t\big) - \P\big(\|G^*\|\op \leq t\big|X\big)\bigg|\\[0.3cm]
& \ \leq \ \sup_{\mathcal{B}\in \mathscr{B}}\bigg|\P\big(G\in \mathcal{B}\big) - \P\big(G^*\in\mathcal{B}\big|X\big)\bigg|\\[0.3cm]
& \ \leq \ \ts\frac{1}{\sqrt 2}\sqrt{d_{\text{KL}}\big(\mathcal{L}(G^*|X)\big\|\mathcal{L}(G)\big)}.
\end{split}
\end{equation*}
Since the distributions $\mathcal{L}(G^*|X)$ and $\mathcal{L}(G)$ are centered and Gaussian, the KL-divergence between them can be calculated explicitly~\citep[][p.33]{pardo2018statistical} in terms of the their respective covariance operators $\hat{\mathsf{C}}$ and $\mathsf C$, 
\begin{equation}\label{eqn:KLformula}
\begin{split}
d_{\text{KL}}\big(\mathcal{L}(G^*|X)\big\|\mathcal{L}(G)\big)
& \ = \ \ts\frac{1}{2}\Big(\tr(\mathsf{C}^{-1}\hat{\mathsf{C}}) \, - \, \tr(\mathsf{I}) \, - \, \log \det(\mathsf{C}^{-1}\hat{\mathsf{C}})\Big).
 \end{split}
\end{equation}
To modify the formula for the KL divergence, recall the formulas for $\mathsf{C}$ and $\hat{\mathsf{C}}$ given in \eqref{eqn:Cinvertibility} and~\eqref{eqn:hatCformula}, and note that the following algebraic relations hold
\begin{equation*}
\begin{split}
\tr(\mathsf{C}^{-1}\hat{\mathsf{C}}) & \ = \ \tr(\hat{\mathsf{D}}\,+\,\mathsf{I})\\[0.2cm]
 \log\det(\mathsf{C}^{-1}\hat{\mathsf{C}} ) & \ = \ \log \det(\hat{\mathsf{D}}\,+\,\mathsf{I}).
 \end{split}
 \end{equation*}
Applying these relations to~\eqref{eqn:KLformula} gives
\begin{equation*}
\begin{split}
d_{\text{KL}}\big(\mathcal{L}(G^*|X)\big\|\mathcal{L}(G)\big)
 & \ = \ \displaystyle\sum_{j=1}^{k(k+1)/2} \ts\frac{1}{2}\Big(\lambda_j(\hat{{\mathsf{D}}})-\log(\lambda_j(\hat{{\mathsf{D}}})+1)\Big).
 \end{split}
\end{equation*}
Noting that $|\lambda_j(\hat{{\mathsf{D}}})|\leq \|\hat{{\mathsf{D}}}\|\op\leq 1/2$, the claim~\eqref{eqn:ivclaim} follows from the basic inequality  $x-\log(x+1)\leq x^2$ for $|x|\leq 1/2$.

Finally, to complete the proof of the proposition, it remains to develop a high-probability upper bound on $\|\hat{\mathsf{D}}\|\op$. By Lemma~\ref{lem:Chatinverse} and Assumption~\ref{A:model}(c), there is a constant $c>0$ not depending on $n$ such that the following bounds hold with probability at least $1-c/n$,
\begin{equation}\label{eqn:tensevenlast}
\begin{split}
 \|\hat{\mathsf{D}}\|\op & \ = \ \big\|\J_k^{-1/2}\big(\,\hat{\!\J}_k \, - \, \J_k\big) \, \J_k^{-1/2} \big\|\op\\[0.2cm]
   & \ \leq \ \ts\frac{1}{\lambda_{\min}(\J_k)}\,\big\|\,\hat{\!\J}_k \, - \, \J_k\big\|\op\\[0.2cm]
 & \ \leq \ \frac{  c\,k \log(n)^{2\nu+1/2}}{n^{1/2}}.
\end{split}
\end{equation}
Applying the bound $k\lesssim n^{\frac{1}{2\beta+4}}$ to~\eqref{eqn:tensevenlast} and~\eqref{eqn:ivclaim} implies the stated result. (It is worth noting that~\eqref{eqn:tensevenlast} actually implies a bound on $\hat{\III}_n$ that is somewhat stronger than the one stated~\eqref{eqn:propIVstatement}. Nevertheless,~\eqref{eqn:propIVstatement} is easier to connect to the rest of the proof of Theorem~\ref{thm:main} and is all that is needed.)\qed

\section{The term $\hat{\I}_n$: Dimension reduction}\label{sec:hatI}
\begin{proposition}
Suppose the conditions of Theorem~\ref{thm:main} hold. Then, there is a constant $c>0$ not depending on $n$ such that the event
\begin{equation*}
\hat{\I}_n \ \leq \ c\,n^{-\frac{\beta-1/2}{2\beta+4+\e}}
\end{equation*}
holds with probability at least $1-c/n$.
\end{proposition}
\proof As a temporary shorthand, let $V^*=\sqrt n \|\hat\Sigma^*-\hat\Sigma\|\op$ and $U^*=\sqrt n \|\hat\Sigma_k^*-\hat\Sigma_k\|\op$, and note that~\eqref{eqn:coupleanti} implies the following almost-sure bound on $\hat{\I}_n$ for any $\delta>0$,
\begin{equation}\label{eqn:reversecoupleanti2}
\hat{\I}_n \ \leq \  \P\big(|V^*-U^*|\geq \delta\,\big|X\big)  \ + \ \sup_{t\in\R}\P\big(|U^*-t|\leq \delta\,|X\big).
\end{equation}
The coupling probability on the right side is handled by Lemma~\ref{lem:covcouplestar}, which shows that there is a constant $c>0$ not depending on $n$, such that if $\delta=cn^{-\frac{\beta-1/2}{2\beta+4}}\log(n)^{\nu+1/2}$, then the event
\begin{equation}\label{eqn:couplesimple}
\P\big(|V^*-U^*|\geq \delta\,\big|X\big) \ \leq \ \ts\frac{c}{n}
\end{equation}
holds with probability at least $1-c/n$.

With regard to the anti-concentration probability on the right side of~\eqref{eqn:reversecoupleanti2}, the reasoning in the proof of Proposition~\ref{prop:I} leading up to~\eqref{eqn:reversecoupleanti} may be re-used to give
\begin{equation*}
  \sup_{t\in\R}\P\big(|U^*-t|\leq \delta \,\big| X\big)
 \ \leq \ \sup_{t\in\R}\P\big(|\|G\|_{_{\scaleto{\!\triangle}{4pt}}}-t|\leq \delta\big)  \ + \ 2\dK\big((\mathcal{L}(U^*|X),\mathcal{L}(\|G\|_{_{\scaleto{\!\triangle}{4pt}}})\big).
\end{equation*}
For the stated choice of $\delta$, Lemma~\ref{lem:Gaussiananti} shows that $ \sup_{t\in\R}\P\big(|\|G\|_{_{\scaleto{\!\triangle}{4pt}}}-t|\leq \delta\big)$ is at most of order $n^{-\frac{\beta-1/2}{2\beta+4+\e}}$.
Also, the triangle inequality for $\dK$ gives
\begin{equation*}
\dK\big((\mathcal{L}(U^*|X),\mathcal{L}(\|G\|_{_{\scaleto{\!\triangle}{4pt}}})\big) \
\leq \ \hat{\II}_n+\hat{\III}_n+\dK(\|G\|\op,\|G\|_{_{\scaleto{\!\triangle}{4pt}}}),
\end{equation*}
and the three terms in this bound are respectively handled in Propositions~\ref{prop:optriI},~\ref{prop:IVhat}, and~\ref{prop:optri}. Collectively, these results show that the sum of terms is at most $cn^{-\frac{\beta-1/2}{2\beta+4+\e}}$ with probability at least $1-c/n$, for some constant $c>0$ not depending on $n$, which completes the proof.\qed

\begin{lemma}\label{lem:covcouplestar}
Suppose the conditions of Theorem~\ref{thm:main} hold. Then, there is a constant $c>0$ not depending on $n$ such that the following event holds with probability at least $1-c/n$,
\begin{equation}\label{eqn:covcouplestar}
\P\bigg( \Big|\sqrt n \|\hat\Sigma^*-\hat \Sigma\|\op- \sqrt n\|\hat\Sigma_k^*-\hat\Sigma_k\|\op\Big| \ \geq  \ c\, n^{-\frac{\beta-1/2}{2\beta+4}}\log(n)^{\nu+1/2} \ \bigg| X\bigg) \ \leq \  \ts\frac{c}{n}.
\end{equation}
\end{lemma}

\proof
For any $t>0$, let the operator ${\mathsf{A}}(t):\mathbb{H}\to \mathbb{H}\oplus\mathbb{H}$ be as defined in the proof of Lemma~\ref{lem:covcouple}, and let $\xi_1^*(t)={\mathsf{A}}(t)X_1^*$.
 The proof of that lemma can be repeated up to~\eqref{eqn:bigLqbound} to establish the following almost-sure bound with $q=\log(n)\vee 3$,
\begin{equation}\label{eqn:bigLqboundstar}
\begin{split}
\Bigg(\E&\bigg| \|\hat\Sigma^*-\hat\Sigma\|\op  - \|\hat\Sigma_k^*-\hat\Sigma_k\|\op\bigg|^q \, \Bigg|X\Bigg)^{1/q}\\[0.2cm]
& \ \leq \  c\bigg( \ts\frac{\sqrt q}{\sqrt n} \cdot \|\E(\xi_1^*(t)\otimes \xi_1^*(t)|X)\|\op^{1/2}\cdot \big(\E(\|\xi_1^*(t)\|^{2q}|X)\big)^{\frac{1}{2q}}\bigg) \bigvee  \bigg(\ts\frac{q}{n}\cdot\big(\E(\|\xi_1^*(t)\|^{2q}|X)\big)^{\frac{1}{q}}\bigg),
\end{split}
\end{equation}
for some absolute constant $c>0$.
%
The operator norm on the right side satisfies
\begin{equation}\label{eqn:xistartensorop}
\begin{split}
\|\E(\xi_1^*(t)\otimes \xi_1^*(t)|X)\|\op & \ = \ \Big\|\ts\frac{1}{n}\sum_{i=1}^n \xi_i(t)\otimes \xi_i(t)\Big\|\op\\[0.2cm]
& \ \leq \ \Big\|\ts\frac{1}{n}\sum_{i=1}^n \xi_i(t)\otimes \xi_i(t) \ - \ \E\, \xi_i(t)\otimes \xi_i(t) \Big\|\op \ + \ \|\E\,\xi_1(t)\otimes \xi_1(t)\|\op.
\end{split}
\end{equation}
Using~\eqref{eqn:Vop} through~\eqref{eqn:lastdimred} with the choice $t=k^{\beta/2}$, we obtain the following bound with probability at least $1-c/n$,
\begin{equation}\label{eqn:V1V1starop}
\begin{split}
\|\E(\xi_1^*(t)\otimes \xi_1^*(t)|X)\|\op 
& \ \leq \ c\,k^{-\beta}.
\end{split}
\end{equation}
Next, to deal with $\E(\|\xi_1^*(t)\|^{2q}|X)^{\frac{1}{2q}}$, we may use~\eqref{eqn:V1norm} to obtain
\begin{equation}\label{eqn:V1star}
\begin{split}
\Big\|\big(\E(\|\xi_1^*(t)\|^{2q}|X)\big)^{\frac{1}{2q}}\Big\|_{L^{2q}}
& \ = \ \Big(\E \|\xi_1(t)\|^{2q}\Big)^{\frac{1}{2q}}\\
& \ \lesssim \ q^{\nu}\, k^{-\beta/2+1/2}.
\end{split}
\end{equation}
So, Chebyshev's inequality implies that the bound
\begin{equation}\label{eqn:xistarL2q}
 (\E(\|\xi_1^*(t)\|^{2q}|X))^{\frac{1}{2q}} \ \leq \ c\,q^{\nu}\, k^{-\beta/2+1/2}
\end{equation}
holds with probability at least $1-c/n$. Applying~\eqref{eqn:V1V1starop} and~\eqref{eqn:xistarL2q} to~\eqref{eqn:bigLqboundstar} shows that the bound
\small
\begin{equation}\label{eqn:Ihatalmostdone}
\Bigg(\E\bigg| \|\hat\Sigma^*-\hat\Sigma\|\op  - \|\hat\Sigma_k^*-\hat\Sigma_k\|\op\bigg|^q \, \Bigg|X\Bigg)^{1/q}\\[0.2cm]
 \ \leq \  \ \frac{c\,q^{\nu+1/2}\, k^{-\beta+1/2}}{\sqrt n}
\end{equation}
\normalsize
holds with probability at least $1-c/n$. Finally, note that we may assume that $k<\dim(\mathbb{H})$, for otherwise the random variables $ \|\hat\Sigma^*-\hat \Sigma\|\op$ and $\|\hat\Sigma_k^*-\hat\Sigma_k\|\op$ are equal. Consequently, the definition of $k$ implies that $k\geq n^{\frac{1}{2\beta+4}}$. When this lower bound on $k$ is applied to the right side of~\eqref{eqn:Ihatalmostdone}, we can obtain the stated result by using a conditional version of Chebyshev's inequality and rescaling with $\sqrt n$.\qed

\section{Background results}\label{sec:background}

To clarify a bit of terminology in the following result, we say that a random matrix $M\in\mathbb{S}^{d\times d}$ is isotropic if $\cov(\llangle A, M\rrangle, \llangle B, M\rrangle)=\llangle A, B\rrangle$ holds for all $A,B\in\mathbb{S}^{d\times d}$.
\begin{lemma}[\cite{Fang:Koike:Balls}, Theorem 2.1]\label{lem:fangkoike}
Let $M_1,\dots,M_n$ be i.i.d.~random matrices in $\mathbb{S}^{d\times d}$ that are centered and isotropic.  In addition, let $\gamma$ denote the standard Gaussian distribution on $\mathbb{S}^{d\times d}$, and let $\mathscr{C}$ denote the class of all convex Borel subsets of $\mathbb{S}^{d\times d}$.
 Then, there is an absolute constant $c>0$ such that
\begin{equation}\label{eqn:Fang:Koike:matrix}
 \sup_{\mathcal{C}\in\mathscr{C}} \bigg| \P\Big(\ts\frac{1}{\sqrt{n}}(M_1+\dots+M_n)\in \mathcal{C}\Big) \ - \ \gamma(\mathcal{C})\bigg| \ \leq  \ \displaystyle \ts\frac{c}{n^{1/2}} \cdot \Big(\ts\frac{d(d+1)}{2}\Big)^{1/4}\cdot (\E\|M_1\|_F^4)^{1/2} \cdot \big(|\log(\E\|M_1\|_F^4/n)|\vee 1\big).
\end{equation}
\end{lemma}
\normalsize
\noindent\textbf{Remark.} In the paper~\cite{Fang:Koike:Balls}, this result is stated in terms of i.i.d.~random vectors $V_1,\dots,V_n\in\R^d$ that are centered and isotropic. Letting $Z\in\R^d$ denote a standard Gaussian vector, and letting $\mathscr{K}$ denote the class of all convex Borel subsets of $\R^d$, the result states that there is an absolute constant $c>0$ such that

\begin{equation}\label{eqn:Fang:Koike:vector}
 \sup_{\mathcal{K}\in\mathscr{K}} \Big| \P\big(\ts\frac{1}{\sqrt{n}}\sum_{i=1}^n V_i\in \mathcal{K}\big) -\P(Z\in\mathcal{K})\Big| \ \leq  \  \ts \frac{c\,d^{1/4}}{n^{1/2}}\cdot (\E\|V_1\|_2^4)^{1/2}\,\big(|\log(\E\|V_1\|_2^4/n)|\vee 1\big).
\end{equation}
This bound directly implies~\eqref{eqn:Fang:Koike:matrix} due to the fact that $(\R^{d(d+1)/2},\langle\cdot,\cdot\rangle$) and $(\mathbb{S}^{d\times d},\llangle\cdot,\cdot\rrangle)$ are isometric. In particular, for any isometry $\psi:\mathbb{S}^{d\times d}\to\R^{d(d+1)/2}$, it can be checked that if $M_1,\dots,M_n$ are i.i.d.~centered isotropic random matrices in $\mathbb{S}^{d\times d}$, then $\psi(M_1),\dots,\psi(M_n)$ are centered isotropic random vectors in $\R^{d(d+1)/2}$. Furthermore, $\psi$ preserves convexity, and the relation $\|M_1\|_F=\|\psi(M_1)\|_2$ holds almost surely.

\begin{lemma}[\cite{Lopes:Bernoulli}, Theorem 2.2]\label{lem:op}
Let $q\geq 3$, and let $\xi_1,\dots,\xi_n$ be i.i.d.~random elements of a real separable Hilbert space with norm $\|\cdot\|$. Also, define the quantity
\begin{equation*}\label{eqn:newrdef}
r(q) \ = \ \frac{q\,\Big(\E\|\xi_1\|^{2q}\Big)^{\frac{1}{q}}}{\big\|\E\,\xi_1\otimes \xi_1\big\|\op}.
\end{equation*}
 Then, there is an absolute constant $c>0$ such that
 \small
\begin{equation*}\label{eqn:hilbresult}
\begin{split}
\bigg(\E\bigg\|\ts\frac{1}{n}\displaystyle\sum_{i=1}^n \xi_i\otimes \xi_i -\E \xi_i\otimes \xi_i\bigg\|\op^q\bigg)^{1/q} 
& \ \leq  \ c\,\big\|\E \xi_1\otimes \xi_1\big\|\op\Big(\sqrt{\ts\frac{r(q)}{n^{1-3/q}}}\,\vee\,\ts\frac{r(q)}{n^{1-3/q}}\Big).
\end{split}
\end{equation*}
\normalsize
In particular, when $q\geq \log(n)\vee 3$, 
\small
\begin{equation*}
\bigg(\E\bigg\|\ts\frac{1}{n}\displaystyle\sum_{i=1}^n \xi_i\otimes \xi_i -\E \xi_i\otimes \xi_i\bigg\|\op^q\bigg)^{1/q}  
\ \leq \ c\bigg( \ts\frac{\sqrt q}{\sqrt n}\, \|\E\xi_1\otimes \xi_1\|\op^{1/2}\, (\E\|\xi_1\|^{2q})^{\frac{1}{2q}}\bigg) \bigvee  \bigg(\ts\frac{q}{n}\, (\E \|\xi_1\|^{2q})^{\frac{1}{q}}\bigg).
\end{equation*}
\end{lemma}
\normalsize
\noindent The following result is a version of Rosenthal's inequality, which is a consequence of Theorem~1 in~\cite{talagrandrosenthal}.
\begin{lemma}\label{lem:rosenthal}
Fix $q\geq 2$, and let $\xi_1,\dots,\xi_n$ be independent centered random variables. Then, there is an absolute constant $c>0$ such that
	\begin{equation*}
	\bigg\|\tsum_{i=1}^n \xi_i\bigg\|_{L^q} \ \leq \ c\,q \, \bigg(\,\bigg\|\tsum_{i=1}^n \xi_i\bigg\|_{L^2} \!\ts\bigvee\, \displaystyle\Big(\tsum_{i=1}^n \|\xi_i\|_{L^q}^q\Big)^{1/q}\bigg).
	\end{equation*}
\end{lemma}

\bibliography{big_alpha_bib_second_revision}

\end{document}